\input amstex
\documentstyle{amsppt}
%
\catcode`@=11
\redefine\output@{%
  \def\break{\penalty-\@M}\let\par\endgraf
  \ifodd\pageno\global\hoffset=105pt\else\global\hoffset=8pt\fi  
  \shipout\vbox{%
    \ifplain@
      \let\makeheadline\relax \let\makefootline\relax
    \else
      \iffirstpage@ \global\firstpage@false
        \let\rightheadline\frheadline
        \let\leftheadline\flheadline
      \else
        \ifrunheads@ 
        \else \let\makeheadline\relax
        \fi
      \fi
    \fi
    \makeheadline \pagebody \makefootline}%
  \advancepageno \ifnum\outputpenalty>-\@MM\else\dosupereject\fi
}
\def\Beta{\mathchar"0\hexnumber@\rmfam 42}
\catcode`\@=\active
\nopagenumbers
\chardef\textvolna='176

\chardef\bigalpha='013
\def\negskp{\hskip -2pt}

\chardef\degree="5E

\def\blue#1{#1}

\edef\darkred#1{#1}
\catcode`#=11\def\diez{#}\catcode`#=6
\catcode`&=11\catcode`&=4
\catcode`_=11\catcode`_=8
\catcode`\^=11\catcode`\^=7
\catcode`~=11\catcode`~=\active
\def\mycite#1{\cite{\blue{#1}}\immediate\special{ps:
     ShrHPSdict begin /ShrBORDERthickness 0 def}}

\def\mytag#1{%
    \tag#1}
\def\mythetag#1{\thetag{\blue{#1}}\immediate\special{ps:
     ShrHPSdict begin /ShrBORDERthickness 0 def}}
\def\myrefno#1{\no#1}
\def\myhref#1#2{\blue{#2}\immediate\special{ps:
     ShrHPSdict begin /ShrBORDERthickness 0 def}}
\def\myEarXivlink{\myhref{http://arXiv.org}{http:/\negskp/arXiv.org}}

\def\mytheorem#1{\csname proclaim\endcsname{Theorem #1}}
\def\mytheoremwithtitle#1#2{\csname proclaim\endcsname{Theorem #1#2}}
\def\mythetheorem#1{\blue{#1}\immediate\special{ps:
     ShrHPSdict begin /ShrBORDERthickness 0 def}}
\def\mylemma#1{\csname proclaim\endcsname{Lemma #1}}
\def\mylemmawithtitle#1#2{\csname proclaim\endcsname{Lemma #1#2}}

\def\mycorollary#1{\csname proclaim\endcsname{Corollary #1}}

\def\mydefinition#1{\definition{Definition #1}}
\def\mythedefinition#1{\blue{#1}\immediate\special{ps:
     ShrHPSdict begin /ShrBORDERthickness 0 def}}
\def\myconjecture#1{\csname proclaim\endcsname{Conjecture #1}}
\def\myconjecturewithtitle#1#2{\csname proclaim\endcsname{Conjecture #1#2}}

\def\myproblem#1{\csname proclaim\endcsname{Problem #1}}
\def\myproblemwithtitle#1#2{\csname proclaim\endcsname{Problem #1#2}}


\pagewidth{360pt}
\pageheight{606pt}
\topmatter
\title
On Walter Wyss's no perfect cuboid paper.
\endtitle
\rightheadtext{On no perfect cuboid paper.}
\author
Ruslan Sharipov
\endauthor
\address Bashkir State University, 32 Zaki Validi street, 450074 Ufa, Russia
\endaddress
\email
\myhref{mailto:r-sharipov\@mail.ru}{r-sharipov\@mail.ru}
\endemail
\abstract
     The perfect cuboid problem is an old famous unsolved problem in mathematics
concerning the existence or non-existence of a rectangular parallelepiped whose
edges, face diagonals, and space diagonal are of integer lengths. Recently
Walter Wyss has published a paper claiming a solution of this problem. The
purpose of this paper is to check out Walter Wyss's result.
\endabstract
\endtopmatter
\TagsOnRight
\document

\head
1. Introduction.
\endhead
\parshape 15 0pt 360pt 0pt 360pt 180pt 180pt
180pt 180pt 180pt 180pt 180pt 180pt 180pt 180pt 180pt 180pt 180pt 180pt
180pt 180pt 180pt 180pt 180pt 180pt 180pt 180pt 180pt 180pt
0pt 360pt
     Actually Walter Wyss's paper is a series of three papers (three versions) each 
updating the previous one. The version 1 is entitled ``No perfect cuboid'' 
(see \mycite{1}), the version 2 is entitled ``On perfect cuboids'' (see \mycite{2}), 
and the version 3 is again entitled ``No perfect cuboid'' (see \mycite{3}). There is
also the version 4 in ArXiv (see \mycite{4}), which looks pretty the same as the
version 3. \vadjust{\vskip 60pt\hbox to 0pt{\kern 15pt
\includegraphics{wyss_01.eps}\hss}\vskip -60pt}Here we consider the version 
4 of Walter Wyss's paper which is the most recent by now. In this paper Walter Wyss 
considers leaning boxes (slanted cuboids) which are parallelepipeds with four rectangular 
faces and two faces being parallelograms (see Fig\.~1.1). Such cuboids have three edges, 
four face diagonals and two space diagonals. Regular cuboids correspond to the case $\theta=\pi/2$. Slanted cuboids and their equations are used in proving the ``No perfect cuboid'' claim for the rectangular case $\theta=\pi/2$.
\head
2. Perfect and non-perfect parallelograms.
\endhead
     The study of slanted cuboids in \mycite{4} is based on the study of perfect 
parallelograms which was carried out in the other paper \mycite{5} by Walter Wyss 
(see also \mycite{6}).
\mydefinition{2.1} A parallelogram is called perfect if its sides and diagonals are
of integer lengths. 
\enddefinition
     Let's consider the parallelogram $ABFE$ in Fig~1.1. and denote their sides
and diagonals by $u_1=|AE|$, $u_2=|EF|$, $u_3=|AF|$, $u_4=|EB|$. Then, applying cosine
theorem, we get the following equations
$$
\hskip -2em
\aligned
&u_3^2=u_1^2+u_2^2-u_1\,u_2\,\cos(\theta),\\
&u_4^2=u_1^2+u_2^2-u_1\,u_2\,\cos(\pi-\theta),
\endaligned
\mytag{2.1}
$$
Since $\cos(\pi-\theta)=-\cos(\theta)$, from \mythetag{2.1} we derive the equation
$$
\hskip -2em
u_3^2+u_4^2=2\,u_1^2+2\,u_2^2.
\mytag{2.2}
$$
\mytheorem{2.1} Four positive real numbers $u_1$, $u_2$, $u_3$, $u_4$ represent
lengths of sides and diagonals of some parallelogram if and only if they obey the 
quadratic equation \mythetag{2.2} and the following inequalities:
$$
\hskip -2em
|u_1-u_2|<u_3<u_1+u_2. 
\mytag{2.3}
$$
\endproclaim
     Let's denote through $A$, $B$, $E$, $F$ a quadruple of points on a plane
such that $|AE|=|BF|$. Using the notations $u_1=|AE|$, $u_2=|EF|$, $u_3=|AF|$, 
$u_4=|EB|$, we see that the inequalities \mythetag{2.3} mean that the points 
$A$, $E$, $F$ constitute a non-degenerate triangle $AEF$. If $ABFE$ is a
parallelogram, the triangle $AEF$ is non-degenerate, hence the inequalities 
\mythetag{2.3} are fulfilled. As we see above, the equation \mythetag{2.2} in 
this case is also fulfilled. So the necessity in Theorem~\mythetheorem{2.1} is 
established.\par
     Let's proceed to the sufficiency. Squaring the inequalities \mythetag{2.3}, 
we get
$$
\hskip -2em
u_1^2+u_2^2-2\,u_1\,u_2<u_3^2<u_1^2+u_2^2+2\,u_1\,u_2. 
\mytag{2.4}
$$
Then, applying the equation \mythetag{2.2} to \mythetag{2.4}, we derive the
inequalities 
$$
\hskip -2em
u_1^2+u_2^2-2\,u_1\,u_2<u_4^2<u_1^2+u_2^2+2\,u_1\,u_2. 
\mytag{2.5}
$$
Since $u_1$, $u_2$, $u_3$, $u_3$ are positive, the inequalities \mythetag{2.5} 
are equivalent to
$$
\hskip -2em
|u_1-u_2|<u_4<u_1+u_2. 
\mytag{2.6}
$$
Due to the notations $u_1=|AE|$, $u_2=|EF|$, $u_3=|AF|$, $u_4=|EB|$ and since
$|AE|=|BF|$, the inequalities \mythetag{2.6} mean that the points 
$B$, $F$, $E$ constitute another non-degenerate triangle $BFE$ with $|BF|=|AE|$. 
The cosine of the angle at the node $F$ in this triangle is calculated as follows:
$$
\hskip -2em
\cos(\hat F)=\frac{u_4^2-u_1^2-u^2}{2\,u_1\,u_2}. 
\mytag{2.7}
$$
Similarly for the cosine of the angle at the node $E$ in the triangle $AEF$ we have:
$$
\hskip -2em
\cos(\hat E)=\frac{u_3^2-u_1^2-u^2}{2\,u_1\,u_2}. 
\mytag{2.8}
$$
Applying the equation \mythetag{2.2} to \mythetag{2.7} and \mythetag{2.8}, we easily
derive 
$$
\hskip -2em
\cos(\hat F)=-\cos(\hat E). 
\mytag{2.9}
$$
The cosine equality \mythetag{2.9} means that
$$
\hskip -2em
\hat F=\pi-\hat E. 
\mytag{2.10}
$$\par
\parshape 16 0pt 360pt 0pt 360pt 175pt 185pt
175pt 185pt 175pt 185pt 175pt 185pt 175pt 185pt 175pt 185pt 175pt 185pt 
175pt 185pt 175pt 185pt 175pt 185pt 175pt 185pt 175pt 185pt 175pt 185pt  
0pt 360pt
Let's draw a triangle $AEF$ using the values of its sides $u_1=|AE|$, $u_2=|EF|$,
\linebreak 
$u_3=|AF|$ and relying on the inequality \mythetag{2.3}. Since $|BF|=|AE|=u_1$ is fixed, 
on the plane there are exactly two locations of the point $B$ relative to the triangle 
$AEF$ at which the equality \mythetag{2.10} holds. \vadjust{\vskip 120pt\hbox to 0pt{\kern 
15pt \includegraphics{wyss_02.eps}\hss}\vskip -120pt}They are symmetric to each
other with respect to the line $EF$ (see Fig.~2.1). Only for one of these two locations
the points $A$, $B$, $E$, $F$ form a parallelogram. Choosing this location, we find that
there is a parallelogram the lengths of whose sides and diagonals coincide with the
numbers $u_1$, $u_2$, $u_3$, $u_4$ obeying the equation \mythetag{2.2} and the
inequalities \mythetag{2.3}. The sufficiency in Theorem~\mythetheorem{2.1} is also
established.\par
     Now let's recall that the inequalities \mythetag{2.3} are equivalent to the
inequalities \mythetag{2.6} modulo the equation \mythetag{2.2}. Therefore 
Theorem~\mythetheorem{2.1} can be reformulated as follows.
\mytheorem{2.2} Four positive real numbers $u_1$, $u_2$, $u_3$, $u_4$ represent
lengths of sides and diagonals of some parallelogram if and only if they obey the 
quadratic equation \mythetag{2.2} and the following inequalities:
$$
\hskip -2em
|u_1-u_2|<u_4<u_1+u_2. 
\mytag{2.11}
$$
\endproclaim
There are two more statements equivalent to Theorem~2.1. 
\mytheorem{2.3} Four positive real numbers $u_1$, $u_2$, $u_3$, $u_4$ represent
lengths of sides and diagonals of some parallelogram if and only if they obey the 
quadratic equation \mythetag{2.2} and the following inequalities:
$$
\xalignat 2
&\hskip -2em
u_3<u_1+u_2,
&&u_4<u_1+u_2. 
\mytag{2.12}
\endxalignat
$$
\endproclaim
\mytheorem{2.4} Four positive real numbers $u_1$, $u_2$, $u_3$, $u_4$ represent
lengths of sides and diagonals of some parallelogram if and only if they obey the 
quadratic equation \mythetag{2.2} and the following inequalities:
$$
\xalignat 2
&\hskip -2em
|u_1-u_2|<u_3,
&&|u_1-u_2|<u_4. 
\mytag{2.13}
\endxalignat
$$
\endproclaim
     According to Definition~\mythedefinition{2.1}, perfect parallelograms are those
for which $u_1$, $u_2$, $u_3$, $u_4$ are positive integers. 
\mydefinition{2.2} A parallelogram is called rational if its sides and diagonals are
of rational lengths. 
\enddefinition
     Rational parallelograms are equivalent to perfect ones since we can bring the
quotients representing the rational numbers $u_1$, $u_2$, $u_3$, $u_4$ to the common
denominator and then obtain a perfect parallelogram by multiplying $u_1$, $u_2$, $u_3$, 
$u_4$ by this common denominator. 
\head
3. Rational slanted cuboids. 
\endhead
\mydefinition{3.1} A slanted cuboid (leaning box) is called perfect if its edges, its
face diagonals, and its space diagonals are of integer lengths. 
\enddefinition
\mydefinition{3.2} A slanted cuboid (leaning box) is called rational if its edges, its
face diagonals, and its space diagonals are of rational lengths. 
\enddefinition
     In \mycite{4} Walter Wyss considers rational leaning boxes. Rational leaning boxes
are equivalent to perfect ones for the same reasons as in the case of rational and perfect
parallelograms (see above). Following Walter Wyss in \mycite{4}, let's consider the
rational leaning box shown in Fig\.~1.1. The edges $AD$, $BC$, $FG$, $EH$ of this leaning 
box are perpendicular to the face parallelogram $ABFE$. Using an appropriate scaling factor,
we can bring their length to the unity
$$
\hskip -2em
|AD|=|BC|=|FG|=|EH|=1.
\mytag{3.1}
$$
Apart from \mythetag{3.1} we use the notations
$$
\gather
\hskip -2em
\aligned
&|AE|=|BF|=|CG|=|DH|=u_1,\\
&|AB|=|EF|=|DC|=|HG|=u_2,\\ 
&|AH|=|ED|=|BG|=|FC|=v_1,\\
&|FH|=|EG|=|AC|=|BD|=v_2,\\
\endaligned
\mytag{3.2}\\
\hskip -2em
\aligned
&|AF|=|DG|=u_3,\\
&|EB|=|HC|=u_4,\\
&|AG|=|DF|=v_3,\\
&|EC|=|BH|=v_4,\\
\endaligned
\mytag{3.3}\\
\endgather
$$
From \mythetag{3.1}, \mythetag{3.2}, and \mythetag{3.3} we derive the
slanted cuboid equations
$$
\gather
\hskip -2em
\aligned
&1+u_1^2=v_1^2,\\
&1+u_2^2=v_2^2,\\
&1+u_3^2=v_3^2,\\
&1+u_4^2=v_4^2,\\
\endaligned
\mytag{3.4}\\
\hskip -2em
2\,u_1^2+2\,u_2^2=u_3^2+u_4^2.
\mytag{3.5}
\endgather
$$
These equation coincide with the equations \thetag{8}--\thetag{12} in \mycite{4}. 
\par
The equations \mythetag{3.4} are Pythagoras equations for rectangular triangles
with rational sides. They can be solved in a parametric form:
$$
\xalignat 2
&\hskip -2em
u_k=\frac{1-s_k^2}{2\,s_k},
&&v_k=\frac{1+s_k^2}{2\,s_k}.
\mytag{3.6}
\endxalignat
$$
Here $k=1,\,\ldots,\,4$ and $s_1$, $s_2$, $s_3$, $s_4$ are rational numbers
obeying the inequalities
$$
0<s_k<1\text{, \ for \ }k=1,\,\ldots,\,4.
\mytag{3.7}
$$\par
     Let's denote through $\psi_1$, $\psi_2$, $\psi_3$, $\psi_4$ the angles 
opposite to the sides of the unit length in the rectangular triangles associoated
with the equations \mythetag{3.4}. Then
$$
\xalignat 2
&\hskip -2em
\sin(\psi_k)=\frac{1}{v_k},
&&\cos(\psi_k)=\frac{u_k}{v_k},
\mytag{3.8}
\endxalignat
$$
where $k=1,\,\ldots,\,4$. Since $u_k$ and $v_k$ are rational numbers, sines and cosines
in \mythetag{3.8} are also rational numbers. 
\mydefinition{3.3}An angle $\psi$ is called a Heron angle if both $\sin(\psi)$ and 
$cos(\psi)$ are rational numbers.
\enddefinition
\mydefinition{3.4}An angle $\psi$ is called an Euler angle if $\tan(\psi)$ is  
a rational number.
\enddefinition
These definitions can be found in Appendix A of the paper \mycite{4}. According to
Definition~\mythedefinition{3.3}, the angles $\psi_1$, $\psi_2$, $\psi_3$, $\psi_4$
in \mythetag{3.8} are Heron angles.\par
     The equation \mythetag{3.5} differs from the equations \mythetag{3.4}. It coincides
with the parallelogram equation \mythetag{2.2}. Substituting \mythetag{3.6} into
\mythetag{3.5} and simplifying, we derive 
$$
\hskip -2em
\gathered
s_4^4\,s_1^2\,s_2^2\,s_3^2+s_3^4\,s_1^2\,s_2^2\,s_4^2-2\,s_2^4\,s_1^2\,s_3^2\,s_4^2
-2\,s_1^4\,s_2^2\,s_3^2\,s_4^2\,+\\
+\,4\,s_1^2\,s_2^2\,s_3^2\,s_4^2-2\,s_2^2\,s_3^2\,s_4^2-2\,s_1^2\,s_3^2\,s_4^2
+s_1^2\,s_2^2\,s_4^2+s_1^2\,s_2^2\,s_3^2=0.
\endgathered
\mytag{3.9}
$$
The parallelogram equation \mythetag{3.5} should be complemented with parallelogram
inequalities. The most simple form of them are given by Theorem~\mythetheorem{2.3}.
Substituting \mythetag{3.6} into \mythetag{2.12} and simplifying, we derive 
$$
\hskip -2em
\gathered
s_1\,s_2^2\,s_3+s_1^2\,s_2\,s_3-s_1\,s_2\,s_3^2+s_1\,s_2-s_2\,s_3-s_1\,s_3<0,\\
s_1\,s_2^2\,s_4+s_1^2\,s_2\,s_4-s_1\,s_2\,s_4^2+s_1\,s_2-s_2\,s_4-s_1\,s_4<0.
\endgathered
\mytag{3.10}
$$
\mytheorem{3.1} Each rational slanted cuboid (leaning box) corresponds to some 
quadruple of rational numbers $s_1$, $s_2$, $s_3$, $s_4$ obeying the polynomial 
equation  \mythetag{3.9} and the polynomial inequalities \mythetag{3.7} and 
\mythetag{3.10}.
\endproclaim
\head
4. Parallelogram parametrization. 
\endhead
     Let's return back to the parallelogram $ABFE$ in Fig\.~1.1. Its sides and diagonals 
are rational numbers $|AE|=u_1$, $|EF|=u_2$, $|AF|=u_3$, $|EB|=u_4$. For such a 
parallelogram Walter Wyss introduces two parameters:
$$
\gather
\hskip -2em
m=\frac{2\,u_2+u_3-u_4}{2\,u_1+u_3+u_4},
\mytag{4.1}\\
\vspace{1ex}
\hskip -2em
n=\frac{2\,u_2-u_3+u_4}{2\,u_1+u_3+u_4}
\mytag{4.2}
\endgather
$$
(see (D.8) and (D.9) in Appendix D of \mycite{4}). Both parameters range within
$$
\xalignat 2
&\hskip -2em
0<m<1, 
&&0<n<1,
\mytag{4.3}
\endxalignat
$$
provided the parallelogram equation \mythetag{3.5} and the parallelogram
inequalities \mythetag{2.3}, \mythetag{2.11}, \mythetag{2.12}, \mythetag{2.13}
are fulfilled.\par
     Indeed, if $m\leqslant 0$, then $u_4\geqslant 2\,u_2+u_3$. Combining this
inequality with the inequality $u_4<u_1+u_2$ from \mythetag{2.11}, we derive 
the inequality $u_3<u_2-u_1$ which contradicts the inequality $|u_1-u_2|<u_3$ from \mythetag{2.3}. If $m\geqslant 1$, then we have $2\,u_2+u_3-u_4\geqslant 2\,u_1+u_3+u_4$. 
This inequality reduces to $u_4\leqslant u_2-u_1$ which contradicts the inequality $|u_1-u_2|<u_4$ from \mythetag{2.11}. Thus, the inequalities for $m$ in \mythetag{4.3} 
are proved. The inequalities for $n$ in \mythetag{4.3} can be proved similarly.\par
     If $u_1\geqslant u_2$ and $u_3\to u_1-u_2$, then from \mythetag{3.5} we derive 
$u_4\to u_1+u_2$. Under these conditions $m\to 0$. Conversely, if $u_2\geqslant u_1$
and $u_3\to u_1+u_2$, then from \mythetag{3.5} we derive $u_4\to u_2-u_1$. Under 
these conditions $m\to 1$. This means that all values from the range $0<m<1$ are
taken by the expression \mythetag{4.1}. Similarly one can prove that all values 
from the range $0<n<1$ are taken by the expression \mythetag{4.2}.\par
     Now let's combine \mythetag{4.1} with the equation \mythetag{3.5} and consider
$$
\hskip -2em
\left\{
\aligned
&2\,u_1^2+2\,u_2^2=u_3^2+u_4^2,\\
&\frac{2\,u_2+u_3-u_4}{2\,u_1+u_3+u_4}=m
\endaligned
\right.
\mytag{4.4}
$$ 
as a system of two equations with $u_3$ and $u_4$ treated as unknowns. Resolving 
the equations \mythetag{4.4} with respect to $u_3$ and $u_4$, we get
$$
\hskip -2em
\aligned
&u_3=\frac{2\,m-m^2+1}{m^2+1}\,u_1+\frac{2\,m+m^2-1}{m^2+1}\,u_2,\\
\vspace{1ex}
&u_4=\frac{1-m^2-2\,m}{m^2+1}\,u_1+\frac{2\,m-m^2+1}{m^2+1}\,u_2.
\endaligned
\mytag{4.5}
$$
The formulas \mythetag{4.5} coincide with the formulas (D.19) and (D.20) 
in Appendix D of \mycite{4}. They can be derived with the use of the 
following Maple\footnotemark\ code:
\footnotetext{\ Maple is a trademark of Waterloo Maple Inc.}
\medskip
\parshape 1 20pt 340pt
\noindent
\darkred{{\tt restart;\newline
Eq\_0:=2*u1\^{}2+2*u2\^{}2-u3\^{}2-u4\^{}2=0:\newline
Eq\_m:=m=(2*u2+u3-u4)/(2*u1+u3+u4):\newline
sss:=solve(\{Eq\_0,Eq\_m\},\{u3,u4\}):\newline
assign(sss):\newline
u3:=collect(u3,[u1,u2]);\newline
u4:=collect(u4,[u1,u2]);
}}
\medskip
The formulas \mythetag{4.5} are understood as a rational parametric solution
of the parallelogram equation \mythetag{3.5} with three parameters
$$
\xalignat 3
&u_1>0,
&&u_2>0,
&&0<m<1.
\endxalignat
$$
\par
     Another parametric solution of the equation \mythetag{3.5} is obtained 
with the use of the formula \mythetag{4.2}. Combining it with \mythetag{3.5},
we write
$$
\hskip -2em
\left\{
\aligned
&2\,u_1^2+2\,u_2^2=u_3^2+u_4^2,\\
&\frac{2\,u_2-u_3+u_4}{2\,u_1+u_3+u_4}=n
\endaligned
\right.
\mytag{4.6}
$$ 
and treat \mythetag{4.6} as a system of equations for unknowns $u_3$ and $u_4$.
Resolving the equations \mythetag{4.6} with respect to $u_3$ and $u_4$, we get
$$
\hskip -2em
\aligned
&u_3=\frac{1-2\,n-n^2}{n^2+1}\,u_1+\frac{1+2\,n-n^2}{n^2+1}\,u_2,\\
\vspace{1ex}
&u_4=\frac{1-n^2+2\,n}{n^2+1}\,u_1+\frac{2\,n+n^2-1}{n^2+1}\,u_2.
\endaligned
\mytag{4.7}
$$
The formulas \mythetag{4.7} coincide with the formulas (D.25) and (D.26) 
in Appendix D of \mycite{4}. They can be derived with the use of the 
following Maple code:
\medskip
\parshape 1 20pt 340pt
\noindent
\darkred{{\tt restart;\newline
Eq\_0:=2*u1\^{}2+2*u2\^{}2-u3\^{}2-u4\^{}2=0:\newline
Eq\_n:=n=(2*u2-u3+u4)/(2*u1+u3+u4):\newline
sss:=solve(\{Eq\_0,Eq\_n\},\{u3,u4\}):\newline
assign(sss):\newline
u3:=collect(u3,[u1,u2]);\newline
u4:=collect(u4,[u1,u2]);
}}
\medskip
The formulas \mythetag{4.7} provide a rational parametric solution
of the parallelogram equation \mythetag{3.5} with three parameters
$$
\xalignat 3
&u_1>0,
&&u_2>0,
&&0<n<1.
\endxalignat
$$
\par
     In Appemndix A of his paper \mycite{4} Walter Wyss introduces the term generator
for an angle. Here is the definition of this term. 
\mydefinition{4.1} For an arbitrary angle $\alpha$ its generator $m=m(\alpha)$ is defined 
by the formula $m=\tan(\alpha/2)$. 
\enddefinition
Any angle $=-\pi<\alpha<\pi$ is uniquely defined by its generator. From the 
formulas
$$
\xalignat 3
&\cos(\alpha)=\frac{1-\tan^2(\alpha/2)}{1+\tan^2(\alpha/2)},
&&\sin(\alpha)=\frac{2\,\tan(\alpha/2)}{1+\tan^2(\alpha/2)},
&&\tan(\alpha/2)=\frac{\sin(\alpha)}{1+\cos(\alpha)},
\endxalignat
$$
which are elementary, we can derive the following theorems.
\mytheorem{4.1} An angle $\alpha\neq\pm\pi$ is a Heron angle if and only if its 
generator is a rational number. 
\endproclaim
\mytheorem{4.2} If an angle $\alpha$ is an Euler angle, then $2\,\alpha$, 
$2\,\alpha-\pi$, and $\pi-2\,\alpha$ are Heron angles. 
\endproclaim
    Using the rational parameters $m$ and $n$ from \mythetag{4.1} and \mythetag{4.2}
as generators, for each rational parallelogram Walter Wyss defines two Heron angles $0<\alpha<\pi/2$ and $0<\beta<\pi/2$ such that
$$
\xalignat 2
&\hskip -2em
\tan(\alpha/2)=m,
&&\tan(\beta/2)=n.
\mytag{4.8}
\endxalignat
$$
Then he introduces two Euler angles
$$
\xalignat 2
&\hskip -2em
\sigma=\frac{\alpha+\beta}{2},
&&\delta=\frac{\alpha-\beta}{2}
\mytag{4.9}
\endxalignat
$$
(see (D.29), (D.30), and (D.39) in Appendix D of \mycite{4}). The latter two
angles obey the inequalities 
$$
\xalignat 2
&0<\sigma<\frac{\pi}{2},
&&-\frac{\pi}{4}<\delta<\frac{\pi}{4}
\endxalignat
$$
(see (D.52) and (D.53) in Appendix D of \mycite{4}).\par
The functions $\omega_{+}$ and $\omega_{-}$ from Appendix C of \mycite{4}
are just notations: 
$$
\xalignat 2
&\hskip -2em
\omega_{+}(x)=\cos(x)+\sin(x),
&&\omega_{-}(x)=\cos(x)-\sin(x).
\mytag{4.10}
\endxalignat
$$
Various formulas using these functions in Appendix D of \mycite{4} can be 
verified with the use of a rational parametrization, e\.\,g\. with the
use of \mythetag{4.7}. Substituting \mythetag{4.7} into the formula
\mythetag{4.1}, upon simplifying we get
$$
\hskip -2em
m=\frac{u_2-n\,u_1}{u_1+n\,u_2}.
\mytag{4.11}
$$
The Maple code responsible for this operation is
\medskip
\parshape 1 20pt 340pt
\noindent
\darkred{{\tt m:=(2*u2+u3-u4)/(2*u1+u3+u4): \newline
m:=normal(m);}}
\medskip
\noindent This code continues the above code on page 7. Therefore the 
\darkred{{\tt restart}} instruction is not issued in it.\par
     Now we need to code the values of sine, cosine, and tangent functions. This is 
done according to \mythetag{4.8}, \mythetag{4.9}, \mythetag{4.10}, and \mythetag{4.11}:
\medskip
\parshape 1 20pt 340pt
\noindent
\darkred{{\tt unprotect(sin,cos,tan):\newline
sin:=subsop(3=NULL,eval(sin)):\newline
cos:=subsop(3=NULL,eval(cos)):\newline
tan:=subsop(3=NULL,eval(tan)):\newline
tan(alpha/2):=m:\newline
tan(beta/2):=n:\newline
sin(alpha/2):=m*cos(alpha/2):\newline
sin(beta/2):=n*cos(beta/2):\newline
sin(alpha):=normal(2*m/(1+m\^{}2)):\newline
cos(alpha):=normal((1-m\^{}2)/(1+m\^{}2)):\newline
sin(beta):=2*n/(1+n\^{}2):\newline
cos(beta):=(1-n\^{}2)/(1+n\^{}2):\newline
omega\_plus:=proc(x) sin(x)+cos(x) end proc:\newline
omega\_minus:=proc(x) cos(x)-sin(x) end proc:
}}
\medskip
\noindent This code continues the previous code and therefore, again, the 
\darkred{{\tt restart}} instruction is not issued in it. Upon running this 
code we can proceed to verifying formulas in Appendix D of \mycite{4}.
In the case of (D.31) we use the following code:
\medskip
\parshape 1 20pt 340pt
\noindent
\darkred{{\tt Expr\_1:=omega\_plus(alpha)*u1-omega\_minus(alpha)*u2-u3:\newline
Expr\_2:=omega\_minus(alpha)*u1+omega\_plus(alpha)*u2-u4:\newline
Expr\_1:=normal(Expr\_1):\newline
Expr\_2:=normal(Expr\_2):\newline
Expr\_1,Expr\_2;
}}
\medskip
\noindent
The output of this code should look like \blue{{\tt 0,0}} confirming that both
expressions \darkred{{\tt Expr\_1}} and \darkred{{\tt Expr\_2}} are zero.\par
     The code verifying the formula (D.32) in Appendix D of \mycite{4} looks
very similar to the previous code. In this case we have 
\medskip
\parshape 1 20pt 340pt
\noindent
\darkred{{\tt Expr\_1:=omega\_plus(alpha)*u3+omega\_minus(alpha)*u4-2*u1:\newline
Expr\_2:=-omega\_minus(alpha)*u3+omega\_plus(alpha)*u4-2*u2:}}\newline
\darkred{{\tt Expr\_1:=normal(Expr\_1):\newline
Expr\_2:=normal(Expr\_2):\newline
Expr\_1,Expr\_2;
}}
\medskip
\noindent
with the same output \blue{{\tt 0,0}} confirming that \darkred{{\tt Expr\_1}} 
and \darkred{{\tt Expr\_2}} both are zero.\par
     The next are the formulas (D.33) and (D.34). They are verified by the code
\medskip
\parshape 1 20pt 340pt
\noindent
\darkred{{\tt Expr\_1:=(u1*u3+u2*u4)/(u1\^{}2+u2\^{}2)-omega\_plus(alpha):\newline
Expr\_2:=(u1*u4-u2*u3)/(u1\^{}2+u2\^{}2)-omega\_minus(alpha):\newline
Expr\_1:=normal(Expr\_1):\newline
Expr\_2:=normal(Expr\_2):\newline
Expr\_1,Expr\_2;
}}
\medskip
    The formulas (D.35) and (D.36) in Appendix D of \mycite{4} are similar
to (D.31) and (D.32). They are verified by the following two fragments of code:
\medskip
\parshape 1 20pt 340pt
\noindent
\darkred{{\tt Expr\_1:=omega\_minus(beta)*u1+omega\_plus(beta)*u2-u3:\newline
Expr\_2:=omega\_plus(beta)*u1-omega\_minus(beta)*u2-u4:\newline
Expr\_1:=normal(Expr\_1):\newline
Expr\_2:=normal(Expr\_2):\newline
Expr\_1,Expr\_2\newline
\newline
Expr\_1:=omega\_minus(beta)*u3+omega\_plus(beta)*u4-2*u1:\newline
Expr\_2:=omega\_plus(beta)*u3-omega\_minus(beta)*u4-2*u2:\newline
Expr\_1:=normal(Expr\_1):\newline
Expr\_2:=normal(Expr\_2):\newline
Expr\_1,Expr\_2;
}}
\medskip
    The formulas (D.37) and (D.38) in Appendix D of \mycite{4} are similar
to (D.33) and (D.34). We use the following code to verify them:
\medskip
\parshape 1 20pt 340pt
\noindent
\darkred{{\tt Expr\_1:=(u1*u4+u2*u3)/(u1\^{}2+u2\^{}2)-omega\_plus(beta):\newline
Expr\_2:=(u1*u3-u2*u4)/(u1\^{}2+u2\^{}2)-omega\_minus(beta):\newline
Expr\_1:=normal(Expr\_1):\newline
Expr\_2:=normal(Expr\_2):\newline
Expr\_1,Expr\_2;
}}
\medskip
    The formulas (D.41), (D.42), and (D.43) in Appendix D of \mycite{4} are 
immediate from (D.31) and (D.35). Therefore we omit their verification and proceed
to (D.44), (D.45), (D.46). These formulas are verified by the following code:
\medskip
\parshape 1 20pt 340pt
\noindent
\darkred{{\tt sigma:=(alpha+beta)/2:\newline
delta:=(alpha-beta)/2:\newline
Expr\_1:=u1*sin(sigma)-u2*cos(sigma):\newline
Expr\_2:=(u1*cos(sigma)+u2*sin(sigma))*omega\_plus(delta)-u3:\newline
Expr\_3:=(u1*cos(sigma)+u2*sin(sigma))*omega\_minus(delta)-u4:\newline
Expr\_1:=normal(expand(Expr\_1)):\newline
Expr\_2:=normal(expand(Expr\_2)):\newline
Expr\_3:=normal(expand(Expr\_3)):\newline
Expr\_2:=subs(cos(1/2*alpha)\^{}2=(cos(alpha)+1)/2,\newline
\phantom{1111111} cos(1/2*beta)\^{}2=(cos(beta)+1)/2,Expr\_2):\newline
Expr\_3:=subs(cos(1/2*alpha)\^{}2=(cos(alpha)+1)/2,\newline
\phantom{1111111} cos(1/2*beta)\^{}2=(cos(beta)+1)/2,Expr\_3):
}}\newline
\darkred{{\tt Expr\_2:=normal(Expr\_2):\newline
Expr\_3:=normal(Expr\_3):\newline
Expr\_1,Expr\_2,Expr\_3;
}}
\medskip
    The formulas (D.47), (D.48), (D.49), (D.50), and (D.51) in Appendix D of 
\mycite{4} are immediate from (D.44), (D.45), (D.46). Therefore we omit their 
verification.\par
    For the reader's convenience all of the above Maple code is placed ibto
the ancillary file \darkred{{\tt section\_04.txt}} attached to this submission.
\par
\head
5. Slanted cuboid formulas. 
\endhead
    The parallelogram equation \mythetag{3.5} is written for the parallelogram
$ABFE$ in Fig\.~1.1. Apart from $ABFE$ there are two other parallelograms 
associated with the slanted cuboid $ABCDEFGH$. They are $AEGC$ and $EFCD$. 
The sides and diagonals of the parallelogram $AEGC$ are
$$
\xalignat 4
&\hskip -2em
|AE|=u_1,
&&|EG|=v_2,
&&|AG|=v_3, 
&&|EC|=v_4.
\qquad
\mytag{5.1}
\endxalignat
$$
Due to \mythetag{5.1} the parallelogram equation for the parallelogram $AEGC$
looks like
$$
\hskip -2em
2\,u_1^2+2\,v_2^2=v_3^2+v_4^2.
\mytag{5.2}
$$
The sides and diagonals of the parallelogram $EFCD$ are
$$
\xalignat 4
&\hskip -2em
|EF|=u_2,
&&|FC|=v_1,
&&|FD|=v_3,
&&|EC|=v_4. 
\qquad
\mytag{5.3}
\endxalignat
$$
Due to \mythetag{5.3} the parallelogram equation for the parallelogram $EFCD$
looks like
$$
\hskip -2em
2\,u_2^2+2\,v_1^2=v_3^2+v_4^2.
\mytag{5.4}
$$
Combining \mythetag{5.2} and \mythetag{5.4} with \mythetag{3.5}, we get  
a system of three equations:
$$
\hskip -2em
\aligned
&2\,u_1^2+2\,u_2^2=u_3^2+u_4^2,\\
&2\,u_1^2+2\,v_2^2=v_3^2+v_4^2,\\
&2\,u_2^2+2\,v_1^2=v_3^2+v_4^2.
\endaligned
\mytag{5.5}
$$
The equations \mythetag{5.5} coincide with (16), (17), and (18) in \mycite{4}. 
\par
{\bf Remark}. The equations \mythetag{5.5} follow from the slanted cuboid equations
\mythetag{3.4} and \mythetag{3.5},  but they are not equivalent to \mythetag{3.4} 
and \mythetag{3.5}.\par
     The equations \mythetag{5.5} are parallelogram equations. Applying the formula
\mythetag{4.1} to them, Walter Wyss defines three rational numbers $m$, $m_1$, $m_2$:
$$
\gather
\hskip -2em
m=\frac{2\,u_2+u_3-u_4}{2\,u_1+u_3+u_4},
\mytag{5.6}\\
\vspace{1ex}
\hskip -2em
m_1=\frac{2\,v_2+v_3-v_4}{2\,u_1+v_3+v_4},
\mytag{5.7}\\
\vspace{1ex}
\hskip -2em
m_2=\frac{2\,u_2+v_3-v_4}{2\,v_1+v_3+v_4}.
\mytag{5.8}
\endgather
$$
These numbers obey the inequalities $0<m<1$, $0<m_1<1$, $0<m_2<1$. Therefore 
they generate three Heron angles $\alpha$, $\alpha_1$, $\alpha_2$ such that
$$
\xalignat 3
&\hskip -2em
0<\alpha<\frac{\pi}{2},
&&0<\alpha_1<\frac{\pi}{2},
&&0<\alpha_2<\frac{\pi}{2}.
\mytag{5.9}
\endxalignat 
$$
\par
     Now we proceed to the formulas \mythetag{3.6}. They are coded as follows:
\medskip
\parshape 1 20pt 340pt
\noindent
\darkred{{\tt restart:\newline 
\newline 
v1:=(s1+1/s1)/2: v2:=(s2+1/s2)/2:\newline
v3:=(s3+1/s3)/2: v4:=(s4+1/s4)/2:\newline
u1:=(1/s1-s1)/2: u2:=(1/s2-s2)/2:\newline 
u3:=(1/s3-s3)/2: u4:=(1/s4-s4)/2:\newline  
\newline 
v1:=normal(v1):  v2:=normal(v2):\newline 
v3:=normal(v3):  v4:=normal(v4):\newline 
u1:=normal(u1):  u2:=normal(u2):\newline 
u3:=normal(u3):  u4:=normal(u4):
}}
\medskip
\noindent
The cuboid equations \mythetag{3.4} are verified by substitution:
\medskip
\parshape 1 20pt 340pt
\noindent
\darkred{{\tt Eq\_1:=1+u1\^{}2-v1\^{}2: Eq\_2:=1+u2\^{}2-v2\^{}2:\newline 
Eq\_3:=1+u3\^{}2-v3\^{}2: Eq\_4:=1+u4\^{}2-v4\^{}2:\newline 
\newline 
normal(Eq\_1), normal(Eq\_2), normal(Eq\_3), normal(Eq\_4);
}}
\medskip
\noindent
The expected output is \blue{{\tt 0,0,0,0}}. It indicates that the equations
\mythetag{3.4} are verified by substituting \mythetag{3.6} into them.\par
     The next step is to derive the equation \mythetag{3.9}. This is done by
the following code:
\medskip
\parshape 1 20pt 340pt
\noindent
\darkred{{\tt Eq\_5:=2*u1\^{}2+2*u2\^{}2-u3\^{}2-u4\^{}2:\newline
Eq\_5:=numer(normal(Eq\_5));
}}
\medskip
The formulas \mythetag{5.6}, \mythetag{5.7}, and \mythetag{5.8} are coded as 
follows:
\medskip
\parshape 1 20pt 340pt
\noindent
\darkred{{\tt m:=normal((2*u2+u3-u4)/(2*u1+u3+u4));\newline
m1:=normal((2*v2+v3-v4)/(2*u1+v3+v4));\newline
m2:=normal((2*u2+v3-v4)/(2*v1+v3+v4));
}}
\medskip
\noindent
The rational numbers $m$, $m_1$ and $m_2$ are used as generators for 
three angles $\alpha$, $\alpha_1$, and $\alpha_2$. This fact is coded
as follows:
\medskip
\parshape 1 20pt 340pt
\noindent
\darkred{{\tt unprotect(sin,cos,tan):\newline
sin:=subsop(3=NULL,eval(sin)):\newline
cos:=subsop(3=NULL,eval(cos)):\newline
tan:=subsop(3=NULL,eval(tan)):\newline
tan(alpha/2):=m:\newline
sin(alpha/2):=m*cos(alpha/2):\newline
sin(alpha):=normal(2*m/(1+m\^{}2)):\newline
cos(alpha):=normal((1-m\^{}2)/(1+m\^{}2)):\newline
\newline
tan(alpha1/2):=m1:\newline
sin(alpha1/2):=m1*cos(alpha1/2):
}}
\medskip
\parshape 1 20pt 340pt
\noindent
\darkred{{\tt sin(alpha1):=normal(2*m1/(1+m1\^{}2)):\newline
cos(alpha1):=normal((1-m1\^{}2)/(1+m1\^{}2)):\newline
\newline
tan(alpha2/2):=m2:\newline
sin(alpha2/2):=m2*cos(alpha2/2):\newline
sin(alpha2):=normal(2*m2/(1+m2\^{}2)):\newline
cos(alpha2):=normal((1-m2\^{}2)/(1+m2\^{}2)):
}}
\medskip
Now the special functions $\omega_{+}$ and $\omega_{-}$ and their special
values are to be programmed. This is done by the following code:
\medskip
\parshape 1 20pt 340pt
\noindent
\darkred{{\tt omega\_plus:=proc(x) sin(x)+cos(x) end proc:\newline
omega\_minus:=proc(x) cos(x)-sin(x) end proc:
}}
\medskip
The number $Q$ is expressed through $s_3$ and $s_4$ by means of the formula
$$
\hskip -2em
Q=s_3\,s_4
\mytag{5.10}
$$
(see (28) in \mycite{4}). It is coded by the following line:
\medskip
\parshape 1 20pt 340pt
\noindent
\darkred{{\tt Q:=s3*s4:
}}
\medskip
The functions $H(x)$, $K(x)$, $M(x)$, $N(x)$ are just notations. They are 
defined in Appendix E of \mycite{4}. Here is the code for them:
\medskip
\parshape 1 20pt 340pt
\noindent
\darkred{{\tt H:=proc(x) global omega\_plus,omega\_minus,Q:\newline 
\phantom{aaaaaaa}omega\_minus(x)-Q*omega\_plus(x) end proc:\newline
K:=proc(x) global omega\_plus,omega\_minus,Q:\newline
\phantom{aaaaaaa}omega\_minus(x)+Q*omega\_plus(x) end proc:\newline
M:=proc(x) global omega\_plus,omega\_minus,Q:\newline
\phantom{aaaaaaa}omega\_plus(x)-Q*omega\_minus(x) end proc:\newline
N:=proc(x) global omega\_plus,omega\_minus,Q:\newline
\phantom{aaaaaaa}omega\_plus(x)+Q*omega\_minus(x) end proc:
}}
\medskip
Now we are able to verify the formulas from section 4 in Walter Wyss's paper
\mycite{4}. Let's begin with the formulas (19) and (20): 
\medskip
\parshape 1 20pt 340pt
\noindent
\darkred{{\tt Eq\_19:=2*u1-u3*omega\_plus(alpha)-u4*omega\_minus(alpha):\newline
Eq\_19:=numer(normal(Eq\_19)):\newline
Eq\_19:=rem(Eq\_19,Eq\_5,s1):\newline
\newline
Eq\_20:=2*u2+u3*omega\_minus(alpha)-u4*omega\_plus(alpha):\newline
Eq\_20:=numer(normal(Eq\_20)):\newline
Eq\_20:=rem(Eq\_20,Eq\_5,s1):\newline
\newline
Eq\_19,Eq\_20;
}}
\medskip
\noindent
The expected output of this code is \blue{{\tt 0,0}}. The \darkred{{\tt rem}}
operator used in this code means that the formulas (19) and (20) hold modulo
the equation \mythetag{3.9}. The same is true for all other formulas in
section 4 of the paper \mycite{4}.\par
     The code below verifies the formulas (21), (22), (23), (24):
\medskip
\parshape 1 20pt 340pt
\noindent
\darkred{{\tt Eq\_21:=2*u1-v3*omega\_plus(alpha1)-v4*omega\_minus(alpha1):\newline
Eq\_21:=numer(normal(Eq\_21)):
}}
\medskip
\parshape 1 20pt 340pt
\noindent
\darkred{{\tt Eq\_21:=rem(Eq\_21,Eq\_5,s1):\newline
\newline
Eq\_22:=2*v2+v3*omega\_minus(alpha1)-v4*omega\_plus(alpha1):\newline
Eq\_22:=numer(normal(Eq\_22)):\newline
Eq\_22:=rem(Eq\_22,Eq\_5,s1):\newline
\newline
Eq\_23:=2*v1-v3*omega\_plus(alpha2)-v4*omega\_minus(alpha2):\newline
Eq\_23:=numer(normal(Eq\_23)):\newline
Eq\_23:=rem(Eq\_23,Eq\_5,s1):\newline
\newline
Eq\_24:=2*u2+v3*omega\_minus(alpha2)-v4*omega\_plus(alpha2):\newline
Eq\_24:=numer(normal(Eq\_24)):\newline
Eq\_24:=rem(Eq\_24,Eq\_5,s1):\newline
\newline
Eq\_21,Eq\_22,Eq\_23,Eq\_24;
}}
\medskip
\noindent
The next are the formulas (29) through (34). They are verified as follows:
\medskip
\parshape 1 20pt 340pt
\noindent
\darkred{{\tt Eq\_29:=4*Q*u1-s4*M(alpha)-s3*H(alpha):\newline
Eq\_29:=numer(normal(Eq\_29)):\newline
Eq\_29:=rem(Eq\_29,Eq\_5,s1):\newline
\newline
Eq\_30:=4*Q*u2+s4*K(alpha)-s3*N(alpha):\newline
Eq\_30:=numer(normal(Eq\_30)):\newline
Eq\_30:=rem(Eq\_30,Eq\_5,s1):\newline
\newline
Eq\_31:=4*Q*u1-s4*N(alpha1)-s3*K(alpha1):\newline
Eq\_31:=numer(normal(Eq\_31)):\newline
Eq\_31:=rem(Eq\_31,Eq\_5,s1):\newline
\newline
Eq\_32:=4*Q*v2+s4*H(alpha1)-s3*M(alpha1):\newline
Eq\_32:=numer(normal(Eq\_32)):\newline
Eq\_32:=rem(Eq\_32,Eq\_5,s1):\newline
\newline
Eq\_33:=4*Q*v1-s4*N(alpha2)-s3*K(alpha2):\newline
Eq\_33:=numer(normal(Eq\_33)):\newline
Eq\_33:=rem(Eq\_33,Eq\_5,s1):\newline
\newline
Eq\_34:=4*Q*u2+s4*H(alpha2)-s3*M(alpha2):\newline
Eq\_34:=numer(normal(Eq\_34)):\newline
Eq\_34:=rem(Eq\_34,Eq\_5,s1):\newline
\newline
Eq\_29,Eq\_30,Eq\_31,Eq\_32,Eq\_33,Eq\_34;
}}
\medskip
     Though the equations (35), (36), (37), (38) are derived from the previous
ones, they can be verified in a straightforward manner:
\medskip
\parshape 1 20pt 340pt
\noindent
\darkred{{\tt Eq\_35:=s4*(M(alpha)-N(alpha1))+s3*(H(alpha)-K(alpha1)):\newline
Eq\_35:=numer(normal(Eq\_35)):\newline
Eq\_35:=rem(Eq\_35,Eq\_5,s1):
}}
\medskip
\parshape 1 20pt 340pt
\noindent
\darkred{{\tt Eq\_36:=-8*Q*u1+s4*(M(alpha)+N(alpha1))+s3*(H(alpha)+K(alpha1)):\newline
Eq\_36:=numer(normal(Eq\_36)):\newline
Eq\_36:=rem(Eq\_36,Eq\_5,s1):\newline
\newline
Eq\_37:=-4*Q*s2+s4*(K(alpha)-H(alpha1))-s3*(N(alpha)-M(alpha1)):\newline
Eq\_37:=numer(normal(Eq\_37)):\newline
Eq\_37:=rem(Eq\_37,Eq\_5,s1):\newline
\newline
Eq\_38:=-4*Q/s2-s4*(K(alpha)+H(alpha1))+s3*(N(alpha)+M(alpha1)):\newline
Eq\_38:=numer(normal(Eq\_38)):\newline
Eq\_38:=rem(Eq\_38,Eq\_5,s1):\newline
\newline
Eq\_35,Eq\_36,Eq\_37,Eq\_38;
}}
\medskip
     Let's recall the formulas \mythetag{4.10}. They can be rewritten as follows:  
$$
\xalignat 2
&\hskip -2em
\omega_{+}(x)=\sqrt{2}\,\cos\Bigl(\frac{\pi}{4}-x\Bigr),
&&\omega_{-}(x)=\sqrt{2}\,\sin\Bigl(\frac{\pi}{4}-x\Bigr).
\mytag{5.11}
\endxalignat
$$
The formulas \mythetag{5.11} are verified by means of the following code:
\medskip
\parshape 1 20pt 340pt
\noindent
\darkred{{\tt omega\_plus(x)-expand(sqrt(2)*cos(Pi/4-x)),\newline
omega\_minus(x)-expand(sqrt(2)*sin(Pi/4-x));
}}
\medskip
\noindent
On page 4 of his paper \mycite{4} Walter Wyss presents the formulas
$$
\xalignat 2
&\hskip -2em
\omega_{+}(\sigma_1)=\sqrt{2}\,\cos\psi,
&&\omega_{-}(\sigma_1)=\sqrt{2}\,\sin\psi,
\mytag{5.12}
\endxalignat
$$
where $2\,\sigma_1=\alpha+\alpha_1$. Comparing \mythetag{5.12} with 
\mythetag{5.11}, we conclude
$$
\hskip -2em
\psi=\frac{\pi}{4}-\sigma_1=\frac{\pi}{4}-\frac{\alpha+\alpha_1}{2}.
\mytag{5.13}
$$
From \mythetag{5.13} one easily derives
$$
\hskip -2em
\alpha+\psi=\frac{\pi}{4}+\frac{\alpha-\alpha_1}{2}=\frac{\pi}{4}+\delta_1,
\mytag{5.14}
$$
where $2\,\delta_1=\alpha-\alpha_1$. Substituting $x=\alpha+\psi$ into
\mythetag{5.11} and using \mythetag{5.14}, we get
$$
\xalignat 2
&\hskip -2em
\omega_{+}(\alpha+\psi)=\sqrt{2}\,\cos\delta_1,
&&\omega_{-}(\alpha+\psi)=-\sqrt{2}\,\sin\delta_1.
\mytag{5.15}
\endxalignat
$$
The formulas \mythetag{5.15} coincide with the formulas given by Walter Wyss
on page 4 of his paper \mycite{4}. So, the formula \mythetag{5.13} is a key point
for understanding what is $\psi$. This formula is programmed by the following
code:
\medskip
\parshape 1 20pt 340pt
\noindent
\darkred{{\tt psi:=Pi/4-alpha/2-alpha1/2:
}}
\medskip
     Note that $\alpha$, $\alpha_1$, $\alpha_2$ are Heron angles \mythetag{5.9}
generated by rational numbers \mythetag{5.6}, \mythetag{5.7}, \mythetag{5.8}. 
Therefore we have the following formulas
$$
\xalignat 3
&\cos^2\Bigl(\frac{\alpha}{2}\Bigr)=\frac{1}{1+m^2},
&&\cos^2\Bigl(\frac{\alpha_1}{2}\Bigr)=\frac{1}{1+m_1^2},
&&\cos^2\Bigl(\frac{\alpha_2}{2}\Bigr)=\frac{1}{1+m_2^2}.
\qquad\quad
\mytag{5.16}
\endxalignat
$$
Relying on the formulas \mythetag{5.16} we introduce a simplification procedure. 
It is called \darkred{{\tt psi\_phi\_simplify}}. We define it with the following
code:
\medskip
\parshape 1 20pt 340pt
\noindent
\darkred{{\tt psi\_phi\_simplify:=proc(A) local AA: global m,m1,m2:\newline
\phantom{aa}AA:=subs(cos(alpha/2)\^{}2=1/(1+m\^{}2),A):\newline
\phantom{aa}AA:=subs(cos(alpha1/2)\^{}2=1/(1+m1\^{}2),AA):\newline
\phantom{aa}AA:=subs(cos(alpha2/2)\^{}2=1/(1+m2\^{}2),AA):\newline
\phantom{aa}return AA:\newline
end proc:
}}
\medskip
Using this procedure, we can proceed to verifying further formulas from 
Walter Wyss's paper. For the formulas (39), (40), (41), (42) we
apply the following code: 
\medskip
\parshape 1 20pt 340pt
\noindent
\darkred{{\tt Eq\_39:=s3*cos(psi)*H(alpha+psi)-s4*sin(psi)*K(alpha+psi):\newline
Eq\_39:=expand(Eq\_39):\newline
Eq\_39:=psi\_phi\_simplify(Eq\_39):\newline
Eq\_39:=numer(normal(Eq\_39)):\newline
Eq\_39:=rem(Eq\_39,Eq\_5,s1):\newline
\newline
Eq\_40:=-4*Q*u1+s4*cos(psi)*M(alpha+psi)+s3*sin(psi)*N(alpha+psi):\newline
Eq\_40:=expand(Eq\_40):\newline
Eq\_40:=psi\_phi\_simplify(Eq\_40):\newline
Eq\_40:=numer(normal(Eq\_40)):\newline
Eq\_40:=rem(Eq\_40,Eq\_5,s1):\newline
\newline
Eq\_41:=-2*Q*s2+s4*cos(psi)*K(alpha+psi)+s3*sin(psi)*H(alpha+psi):\newline
Eq\_41:=expand(Eq\_41):\newline
Eq\_41:=psi\_phi\_simplify(Eq\_41):\newline
Eq\_41:=numer(normal(Eq\_41)):\newline
Eq\_41:=rem(Eq\_41,Eq\_5,s1):\newline
\newline
Eq\_42:=-2*Q/s2-s4*sin(psi)*M(alpha+psi)+s3*cos(psi)*N(alpha+psi):\newline
Eq\_42:=expand(Eq\_42):\newline
Eq\_42:=psi\_phi\_simplify(Eq\_42):\newline
Eq\_42:=numer(normal(Eq\_42)):\newline
Eq\_42:=rem(Eq\_42,Eq\_5,s1):\newline
\newline
Eq\_39,Eq\_40,Eq\_41,Eq\_42;
}}
\medskip
\noindent
The formulas (43) and (44) are matrix presentations of the formulas
(39), (40), (41), (42). The formulas (45) and (46) are inverse to (43) and (44).
We do not verify them. However, we do verify the formulas (47), (48), (49), (50)
derived from (45) and (46). This is done by the following code:
\medskip
\parshape 1 20pt 340pt
\noindent
\darkred{{\tt Eq\_47:=-K(alpha+psi)+2*s2*s3*cos(psi):\newline
Eq\_47:=expand(Eq\_47):\newline
Eq\_47:=numer(normal(Eq\_47)):\newline
Eq\_47:=rem(Eq\_47,Eq\_5,s1):\newline
\newline
Eq\_48:=-H(alpha+psi)+2*s2*s4*sin(psi):\newline
Eq\_48:=expand(Eq\_48):
}}
\medskip
\parshape 1 20pt 340pt
\noindent
\darkred{{\tt Eq\_48:=numer(normal(Eq\_48)):\newline
Eq\_48:=rem(Eq\_48,Eq\_5,s1):\newline
\newline
Eq\_49:=-M(alpha+psi)+4*u1*s3*cos(psi)-2*s3/s2*sin(psi):\newline
Eq\_49:=expand(Eq\_49):\newline
Eq\_49:=numer(normal(Eq\_49)):\newline
Eq\_49:=rem(Eq\_49,Eq\_5,s1):\newline
\newline
Eq\_50:=-N(alpha+psi)+4*u1*s4*sin(psi)+2*s4/s2*cos(psi):\newline
Eq\_50:=expand(Eq\_50):\newline
Eq\_50:=numer(normal(Eq\_50)):\newline
Eq\_50:=rem(Eq\_50,Eq\_5,s1):\newline
\newline
Eq\_47,Eq\_48,Eq\_49,Eq\_50;
}}
\medskip
    The formulas (51), (52), (53), (54) in \mycite{4} are similar to the 
formulas (35), (36), (37), (38). They are verified by means of the following
code:
\medskip
\parshape 1 20pt 340pt
\noindent
\darkred{{\tt Eq\_51:=s3*(N(alpha)-M(alpha2))-s4*(K(alpha)-H(alpha2)):\newline
Eq\_51:=numer(normal(Eq\_51)):\newline
Eq\_51:=rem(Eq\_51,Eq\_5,s1):\newline
\newline
Eq\_52:=-8*Q*u2+s3*(N(alpha)+M(alpha2))-s4*(K(alpha)+H(alpha2)):\newline
Eq\_52:=numer(normal(Eq\_52)):\newline
Eq\_52:=rem(Eq\_52,Eq\_5,s1):\newline
\newline
Eq\_53:=-4*Q*s1+s4*(N(alpha2)-M(alpha))+s3*(K(alpha2)-H(alpha)):\newline
Eq\_53:=numer(normal(Eq\_53)):\newline
Eq\_53:=rem(Eq\_53,Eq\_5,s1):\newline
\newline
Eq\_54:=-4*Q/s1+s4*(N(alpha2)+M(alpha))+s3*(K(alpha2)+H(alpha)):\newline
Eq\_54:=numer(normal(Eq\_54)):\newline
Eq\_54:=rem(Eq\_54,Eq\_5,s1):\newline
\newline
Eq\_51,Eq\_52,Eq\_53,Eq\_54;
}}
\medskip
     In the next fragment of Walter Wyss's paper \mycite{4} the angle $\phi$
is defined:
$$
\hskip -2em
\phi=\frac{\pi}{4}-\sigma_2=\frac{\pi}{4}-\frac{\alpha+\alpha_2}{2}.
\mytag{5.17}
$$
Here $2\,\sigma_2=\alpha+\alpha_2$. Though the formula \mythetag{5.17} is not
written explicitly, the formulas 
$$
\xalignat 2
&\hskip -2em
\omega_{+}(\sigma_2)=\sqrt{2}\,\cos\phi,
&&\omega_{-}(\sigma_2)=\sqrt{2}\,\sin\phi,
\mytag{5.18}
\endxalignat
$$
compared with \mythetag{5.11} lead to \mythetag{5.17}. Then the following 
formula with $2\,\delta_2=\alpha-\alpha_2$ is derived from \mythetag{5.17}:
$$
\hskip -2em
\alpha+\phi=\frac{\pi}{4}+\frac{\alpha-\alpha_2}{2}=\frac{\pi}{4}+\delta_2,
\mytag{5.19}
$$
Substituting $x=\alpha+\phi$ into
\mythetag{5.11} and using \mythetag{5.19}, we get
$$
\xalignat 2
&\hskip -2em
\omega_{+}(\alpha+\phi)=\sqrt{2}\,\cos\delta_2,
&&\omega_{-}(\alpha+\phi)=-\sqrt{2}\,\sin\delta_2.
\mytag{5.20}
\endxalignat
$$
The formulas \mythetag{5.20} are analogous to \mythetag{5.15}, the formula
\mythetag{5.19} is analogous to \mythetag{5.14}, the formulas
\mythetag{5.18} are analogous to \mythetag{5.12}, and the formula
\mythetag{5.17} is ana\-logous to \mythetag{5.13}. The formula \mythetag{5.17}
is programmed by the following code:
\medskip
\parshape 1 20pt 340pt
\noindent
\darkred{{\tt phi:=Pi/4-alpha/2-alpha2/2:
}}
\medskip
     The angle $\phi$ is used by Walter Wyss in his formulas (55), (56), (57), 
(58). These formulas are verified as follows:
\medskip
\parshape 1 20pt 340pt
\noindent
\darkred{{\tt Eq\_55:=-s4*cos(phi)*K(alpha+phi)-s3*sin(phi)*H(alpha+phi):\newline
Eq\_55:=expand(Eq\_55):\newline
Eq\_55:=psi\_phi\_simplify(Eq\_55):\newline
Eq\_55:=numer(normal(Eq\_55)):\newline
Eq\_55:=rem(Eq\_55,Eq\_5,s1):\newline
\newline
Eq\_56:=-4*Q*u2-s4*sin(phi)*M(alpha+phi)+s3*cos(phi)*N(alpha+phi):\newline
Eq\_56:=expand(Eq\_56):\newline
Eq\_56:=psi\_phi\_simplify(Eq\_56):\newline
Eq\_56:=numer(normal(Eq\_56)):\newline
Eq\_56:=rem(Eq\_56,Eq\_5,s1):\newline
\newline
Eq\_57:=-2*Q*s1+s4*sin(phi)*K(alpha+phi)-s3*cos(phi)*H(alpha+phi):\newline
Eq\_57:=expand(Eq\_57):\newline
Eq\_57:=psi\_phi\_simplify(Eq\_57):\newline
Eq\_57:=numer(normal(Eq\_57)):\newline
Eq\_57:=rem(Eq\_57,Eq\_5,s1):\newline
\newline
Eq\_58:=-2*Q/s1+s4*cos(phi)*M(alpha+phi)+s3*sin(phi)*N(alpha+phi):\newline
Eq\_58:=expand(Eq\_58):\newline
Eq\_58:=psi\_phi\_simplify(Eq\_58):\newline
Eq\_58:=numer(normal(Eq\_58)):\newline
Eq\_58:=rem(Eq\_58,Eq\_5,s1):\newline
\newline
Eq\_55,Eq\_56,Eq\_57,Eq\_58;
}}
\medskip
We omit the formulas (59), (60), (61), (62) just like the formulas (43), (44), (45), 
(46) above and proceed to (63), (64), (65), (66):
\medskip
\parshape 1 20pt 340pt
\noindent
\darkred{{\tt Eq\_63:=-K(alpha+phi)+2*s1*s3*sin(phi):\newline
Eq\_63:=expand(Eq\_63):\newline
Eq\_63:=numer(normal(Eq\_63)):\newline
Eq\_63:=rem(Eq\_63,Eq\_5,s1):\newline
\newline
Eq\_64:=-H(alpha+phi)-2*s1*s4*cos(phi):\newline
Eq\_64:=expand(Eq\_64):\newline
Eq\_64:=numer(normal(Eq\_64)):\newline
Eq\_64:=rem(Eq\_64,Eq\_5,s1):\newline
\newline
Eq\_65:=-M(alpha+phi)-4*u2*s3*sin(phi)+2*s3/s1*cos(phi):\newline
Eq\_65:=expand(Eq\_65):\newline
Eq\_65:=numer(normal(Eq\_65)):
}}
\medskip
\parshape 1 20pt 340pt
\noindent
\darkred{{\tt 
Eq\_65:=rem(Eq\_65,Eq\_5,s1):\newline
\newline
Eq\_66:=-N(alpha+phi)+4*u2*s4*cos(phi)+2*s4/s1*sin(phi):\newline
Eq\_66:=expand(Eq\_66):\newline
Eq\_66:=numer(normal(Eq\_66)):\newline
Eq\_66:=rem(Eq\_66,Eq\_5,s1):\newline
\newline
Eq\_63,Eq\_64,Eq\_65,Eq\_66;
}}
\medskip
     The next are the formulas (67), (68), (69), (70) in Walter Wyss's paper \mycite{4}. 
They are verified by means of the following code:
\medskip
\parshape 1 20pt 340pt
\noindent
\darkred{{\tt Eq\_67:=-4*Q*u2-s4*K(alpha)+s3*N(alpha):\newline
Eq\_67:=numer(normal(Eq\_67)):\newline
Eq\_67:=rem(Eq\_67,Eq\_5,s1):\newline
\newline
Eq\_68:=-4*Q*v2-s4*H(alpha1)+s3*M(alpha1):\newline
Eq\_68:=numer(normal(Eq\_68)):\newline
Eq\_68:=rem(Eq\_68,Eq\_5,s1):\newline
\newline
Eq\_69:=-4*Q*u1+s4*M(alpha)+s3*H(alpha):\newline
Eq\_69:=numer(normal(Eq\_69)):\newline
Eq\_69:=rem(Eq\_69,Eq\_5,s1):\newline
\newline
Eq\_70:=-4*Q*v1+s4*N(alpha2)+s3*K(alpha2):\newline
Eq\_70:=numer(normal(Eq\_70)):\newline
Eq\_70:=rem(Eq\_70,Eq\_5,s1):\newline
\newline
Eq\_67,Eq\_68,Eq\_69,Eq\_70;
}}
\medskip
There are two formulas on page 8 of the paper \mycite{4}. They are not numbered:
$$
\gather
\hskip -2em
-s_4\,K(\alpha)+s_3\,N(\alpha)=-s_4\,H(\alpha_2)+s_3\,M(\alpha_2),
\mytag{5.21}
\\
\hskip -2em
s_4\,M(\alpha)+s_3\,H(\alpha)=s_4\,N(\alpha_1)+s_3\,K(\alpha_1).
\mytag{5.22}
\endgather
$$
Giving them the numbers \mythetag{5.21} and \mythetag{5.22}, we can verify them
as follows:
\medskip
\parshape 1 20pt 340pt
\noindent
\darkred{{\tt Eq\_5\_21:=-s4*K(alpha)+s3*N(alpha)+s4*H(alpha2)-s3*M(alpha2):\newline
Eq\_5\_21:=numer(normal(Eq\_5\_21)):\newline
Eq\_5\_21:=rem(Eq\_5\_21,Eq\_5,s1):\newline
\newline
Eq\_5\_22:=s4*M(alpha)+s3*H(alpha)-s4*N(alpha1)-s3*K(alpha1):\newline
Eq\_5\_22:=numer(normal(Eq\_5\_22)):\newline
Eq\_5\_22:=rem(Eq\_5\_22,Eq\_5,s1):\newline
\newline
Eq\_5\_21,Eq\_5\_22;
}}
\medskip
We do not need to follow the proof of the equations (50) and (66) on page 8
of \mycite{4}. These equations are verified programmatically above. Similarly,
we do not need to follow the proof of the formula (71) on page 9 of this paper. 
\pagebreak This formula is also verified programmatically by means of the following 
code:
\medskip
\parshape 1 20pt 340pt
\noindent
\darkred{{\tt Eq\_71:=s1*s2-tan(phi-psi):\newline
Eq\_71:=expand(Eq\_71):\newline
Eq\_71:=numer(normal(Eq\_71)):\newline
Eq\_71:=rem(Eq\_71,Eq\_5,s1);
}}
\medskip
\noindent
Applying the formulas \mythetag{5.13} and \mythetag{5.17}, we can write the formula
(71) as follows:
$$
\hskip -2em
\tan\Bigl(\frac{\alpha_1-\alpha_2}{2}\Bigr)=s_1\,s_2.
\mytag{5.23}
$$
The formula \mythetag{5.23} can also be verified programmatically:
\medskip
\parshape 1 20pt 340pt
\noindent
\darkred{{\tt Eq\_5\_23:=s1*s2-tan(alpha1/2-alpha2/2):\newline
Eq\_5\_23:=expand(Eq\_5\_23):\newline
Eq\_5\_23:=numer(normal(Eq\_5\_23)):\newline
Eq\_5\_23:=rem(Eq\_5\_23,Eq\_5,s1);
}}
\medskip
     Section 5 of Walter Wyss's paper \mycite{4} is slightly different from 
section 4. Nevertheless, now we proceed to this section and verify some
prerequisite formulas therein. The formulas (75), (76), (77), (78) are
verified as follows:
\medskip
\parshape 1 20pt 340pt
\noindent
\darkred{{\tt Eq\_75:=-K(alpha)+2*s3*(s2*cos(psi)\^{}2\newline
\phantom{aaaaaaa}+sin(psi)*(2*u1*cos(psi)-1/s2*sin(psi))):\newline
Eq\_75:=expand(Eq\_75):\newline
Eq\_75:=psi\_phi\_simplify(Eq\_75):\newline
Eq\_75:=numer(normal(Eq\_75)):\newline
Eq\_75:=rem(Eq\_75,Eq\_5,s1):\newline
\newline
Eq\_76:=-N(alpha)+2*s4*(-s2*sin(psi)\^{}2\newline
\phantom{aaaaaaa}+cos(psi)*(2*u1*sin(psi)+1/s2*cos(psi))):\newline
Eq\_76:=expand(Eq\_76):\newline
Eq\_76:=psi\_phi\_simplify(Eq\_76):\newline
Eq\_76:=numer(normal(Eq\_76)):\newline
Eq\_76:=rem(Eq\_76,Eq\_5,s1):\newline
\newline
Eq\_77:=-H(alpha)+2*s4*(s2*sin(psi)*cos(psi)\newline
\phantom{aaaaaaa}+sin(psi)*(2*u1*sin(psi)+1/s2*cos(psi))):\newline
Eq\_77:=expand(Eq\_77):\newline
Eq\_77:=psi\_phi\_simplify(Eq\_77):\newline
Eq\_77:=numer(normal(Eq\_77)):\newline
Eq\_77:=rem(Eq\_77,Eq\_5,s1):\newline
\newline
Eq\_78:=-M(alpha)+2*s3*(-s2*sin(psi)*cos(psi)\newline
\phantom{aaaaaaa}+cos(psi)*(2*u1*cos(psi)-1/s2*sin(psi))):\newline
Eq\_78:=expand(Eq\_78):\newline
Eq\_78:=psi\_phi\_simplify(Eq\_78):\newline
Eq\_78:=numer(normal(Eq\_78)):\newline
Eq\_78:=rem(Eq\_78,Eq\_5,s1):\newline
\newline
Eq\_75,Eq\_76,Eq\_77,Eq\_78;
}}
\medskip
    The formula (79) coincides with (29), the formula (80) coincides with (30).
The formula (81) is just a notation. Taking into account this notation, the formulas
(82), (83), (84) are verified as follows:
\medskip
\parshape 1 20pt 340pt
\noindent
\darkred{{\tt lambda:=tan(psi):\newline
\newline
unprotect(cot):\newline
cot:=subsop(3=NULL,eval(cot)):\newline
cot(alpha/2):=1/m:\newline
cot(alpha1/2):=1/m1:\newline
cot(alpha2/2):=1/m2:\newline
\newline
Eq\_82:=-omega\_minus(alpha)+lambda*omega\_plus(alpha)\newline
\phantom{aaaaaaa}+s2*s3+lambda*s2*s4:\newline
Eq\_82:=expand(Eq\_82):\newline
Eq\_82:=numer(normal(Eq\_82)):\newline
Eq\_82:=rem(Eq\_82,Eq\_5,s1):\newline
\newline
Eq\_83:=-Q*(omega\_plus(alpha)+lambda*omega\_minus(alpha))\newline
\phantom{aaaaaaa}+s2*s3-lambda*s2*s4:\newline
Eq\_83:=expand(Eq\_83):\newline
Eq\_83:=numer(normal(Eq\_83)):\newline
Eq\_83:=rem(Eq\_83,Eq\_5,s1):\newline
\newline
Eq\_84:=-2*u1*(omega\_minus(alpha)-lambda*omega\_plus(alpha))\newline
\phantom{aaaaaaa}-s4+lambda*s3+s2*(omega\_plus(alpha)\newline
\phantom{aaaaaaa}+lambda*omega\_minus(alpha)):\newline
Eq\_84:=expand(Eq\_84):\newline
Eq\_84:=numer(normal(Eq\_84)):\newline
Eq\_84:=rem(Eq\_84,Eq\_5,s1):\newline
\newline
Eq\_82,Eq\_83,Eq\_84;
}}
\medskip
\noindent
Concluding the above computations, we can confirm that the formulas (19)--(84)
in Walter Wyss's paper are valid.
\par
\head
6. A special solution of the slanted cuboid equation. 
\endhead
    Theorem~\mythetheorem{3.1} provides an exhaustive description of rational
slanted cuboids. They constitute rational points within an open subvariety of 
a three-dimensional real algebraic variety in $\Bbb R^4$. This real algebraic 
variety $\Gamma_3$ is defined by the equation \mythetag{3.9}. Its open subvariety $\Gamma_{3\sssize ++}\subset\Gamma_3$ is outlined by the inequalities 
\mythetag{3.7} and \mythetag{3.10}.\par
     Let's consider the equality (85) in Walter Wyss's paper \mycite{4}. It
is different from all of the previous formulas in this paper. The equality (85)
does not hold identically on $\Gamma_3$. It makes an auxiliary restriction 
thus defining a two-dimensional subvariety $\Gamma_2^1\subset\Gamma_3$. The 
lower index $2$ in $\Gamma_2^1$ indicates the dimension of the subvariety. The upper
index $1$ in $\Gamma_2^1$ says that $\Gamma_2^1$ is not the only two-dimensional 
subvariety of $\Gamma_3$ that will be considered in what follows. Thus, an auxiliary 
restriction is set:
$$
\hskip -2em
\lambda=\tan\psi=0.
\mytag{6.1}
$$ 
Under the restriction \mythetag{6.1} the formulas (82), (83), (84) in \mycite{4}
reduce to the formulas (86), (87), (88) therein. Here are these formulas 
$$
\gather
\hskip -2em
s_2\,s_3=\omega_{-}(\alpha),
\mytag{6.2}\\
\hskip -2em
s_2\,s_3=Q\,\omega_{+}(\alpha),
\mytag{6.3}\\
\hskip -2em
s_2\,\omega_{+}(\alpha)=2\,u_1\,\omega_{-}(\alpha)+s_4.
\mytag{6.4}
\endgather
$$
Following Walter Wyss in \mycite{4}, we multiply both sides of \mythetag{6.4}
by $\omega_{+}(\alpha)$:
$$
\hskip -2em
s_2\,\omega^2_{+}(\alpha)=2\,u_1\,\omega_{-}(\alpha)\,\omega_{+}(\alpha)
+s_4\,\omega_{+}(\alpha).
\mytag{6.5}
$$
Then we recall the formula \mythetag{5.10} for $Q$. Applying \mythetag{5.10} to
\mythetag{6.3}, we get
$$
\hskip -2em
s_2=s_4\,\omega_{+}(\alpha).
\mytag{6.6}
$$
Due to \mythetag{6.6} we can replace the last term in \mythetag{6.5} with $s_2$:
$$
\hskip -2em
s_2\,\omega^2_{+}(\alpha)=2\,u_1\,\omega_{-}(\alpha)\,\omega_{+}(\alpha)+s_2.
\mytag{6.7}
$$
The formula \mythetag{6.7} can be transformed as 
$$
\hskip -2em
s_2\,(\omega^2_{+}(\alpha)-1)=2\,u_1\,\omega_{-}(\alpha)\,\omega_{+}(\alpha).
\mytag{6.8}
$$
Now we recall the formulas \mythetag{4.10}. From these formulas we derive
$$
\gather
\hskip -2em
\omega^2_{+}(\alpha)-1=\sin(2\,\alpha),
\mytag{6.9}\\
\hskip -2em
\omega_{-}(\alpha)\,\omega_{+}(\alpha)=\cos(2\,\alpha).
\mytag{6.10}
\endgather
$$
Applying \mythetag{6.9} and \mythetag{6.10} to \mythetag{6.8}, we obtain
$$
\hskip -2em
s_2\,\sin(2\,\alpha)=2\,u_1\,\cos(2\,\alpha).
\mytag{6.11}
$$
The formula \mythetag{6.11} is equivalent to the formula (91) in \mycite{4}.
\par
     Now let's recall that the angle $\alpha$ in \mythetag{5.9} was introduced
through its generator \mythetag{5.6} (see Definition~\mythedefinition{4.1}),
i\.\,e\. we have the formula
$$
\hskip -2em
\tan(\alpha/2)=m.
\mytag{6.12}
$$ 
From \mythetag{6.12} we derive 
$$
\gather
\hskip -2em
\cos(\alpha)=\frac{1-\tan^2(\alpha/2)}{1+\tan^2(\alpha/2)}
=\frac{1-m^2}{1+m^2},
\mytag{6.13}\\
\vspace{1ex}
\hskip -2em
\sin(\alpha)=\frac{2\,\tan(\alpha/2)}{1+\tan^2(\alpha/2)}
=\frac{2\,m}{1+m^2}.
\mytag{6.14}
\endgather
$$
Then from \mythetag{6.13} and \mythetag{6.14} we derive 
$$
\hskip -2em
\sin(2\,\alpha)=2\,\sin\alpha\,\cos\alpha=\frac{4\,m\,\bigl(1-m^2\bigr)}
{\bigl(1+m^2\bigr)^2}.
\mytag{6.15}
$$
Again from \mythetag{6.13} and \mythetag{6.14} we derive 
$$
\hskip -2em
\cos(2\,\alpha)=\cos^2\!\alpha-\sin^2\!\alpha
=\frac{\bigl(1-m^2\bigr)^2-4\,m^2}
{\bigl(1+m^2\bigr)^2}.
\mytag{6.16}
$$
Now we apply \mythetag{6.15} and \mythetag{6.16} to \mythetag{6.11}. As a
result we get 
$$
\hskip -2em
s_2=u_1\,\frac{\bigl(1-m^2\bigr)^2-4\,m^2}{2\,m\,\bigl(1-m^2\bigr)}
\mytag{6.17}
$$
\par
    Let's denote $s_1=s$. Then the first formula \mythetag{3.6} for $k=1$ is 
written as
$$
\hskip -2em
u_1=\frac{1-s^2}{2\,s},
\mytag{6.18}
$$
Substituting \mythetag{6.18} into \mythetag{6.17}, we derive the formula
$$
\hskip -2em
s_2=\frac{\bigl(1-s^2\bigr)\,\bigl(\bigl(1-m^2\bigr)^2
-4\,m^2\bigr)}{4\,m\,s\,\bigl(1-m^2\bigr)}.
\mytag{6.19}
$$
Applying \mythetag{4.10}, \mythetag{6.13}, and \mythetag{6.14} once more, we 
obtain the formulas
$$
\gather
\hskip -2em
\omega_{-}(\alpha)=\cos(\alpha)-\sin(\alpha)
=\frac{1-m^2-2\,m}{1+m^2},
\mytag{6.20}\\
\vspace{1ex}
\hskip -2em
\omega_{+}(\alpha)=\cos(\alpha)+\sin(\alpha)
=\frac{1-m^2+2\,m}{1+m^2}.
\mytag{6.21}
\endgather
$$
Now we substitute \mythetag{6.19}, and \mythetag{6.20} into \mythetag{6.2} 
and we get 
$$
\hskip -2em
s_3=\frac{4\,m\,s\,\bigl(1-m^2\bigr)}{\bigl(1-s^2\bigr)
\,\bigl(1+m^2\bigr)\,\bigl(1-m^2+2\,m\bigr)},
\mytag{6.22}
$$
Then we substitute \mythetag{6.19}, and \mythetag{6.21} into \mythetag{6.6}.
As a result we get
$$
\hskip -2em
s_4=\frac{\bigl(1-s^2\bigr)\,\bigl(1+m^2\bigr)\,
\bigl(1-m^2-2\,m\bigr)}{4\,m\,s\,\bigl(1-m^2\bigr)}.
\mytag{6.23}
$$
Let's denote through $\theta(s,m)$, $\eta(s,m)$, $\zeta(s,m)$ the right hand sides
of the formulas \mythetag{6.19}, \mythetag{6.22} and \mythetag{6.23} respectively.
The symbol $m$ is linked with the angle $\alpha$ in \mythetag{5.9}. In order unlink
the argument $m$ of the functions $\theta(s,m)$, $\eta(s,m)$, $\zeta(s,m)$ from  
the angle $\alpha$ we replace it with $\mu$. As a result we have
$$
\pagebreak
\gather
\hskip -2em
\theta(s,\mu)=\frac{\bigl(1-s^2\bigr)\,\bigl(\bigl(1-\mu^2\bigr)^2
-4\,\mu^2\bigr)}{4\,\mu\,s\,\bigl(1-\mu^2\bigr)},
\mytag{6.24}\\
\vspace{1ex}
\hskip -2em
\eta(s,\mu)=\frac{4\,\mu\,s\,\bigl(1-\mu^2\bigr)}{\bigl(1-s^2\bigr)
\,\bigl(1+\mu^2\bigr)\,\bigl(1-\mu^2+2\,\mu\bigr)},
\mytag{6.25}\\
\vspace{1ex}
\hskip -2em
\zeta(s,\mu)=\frac{\bigl(1-s^2\bigr)\,\bigl(1+\mu^2\bigr)\,
\bigl(1-\mu^2-2\,\mu\bigr)}{4\,\mu\,s\,\bigl(1-\mu^2\bigr)}.
\mytag{6.26}
\endgather
$$
Using the functions \mythetag{6.24}, \mythetag{6.25}, \mythetag{6.26}, we define
a mapping:
$$
\hskip -2em
\left\{
\aligned
&s_1=s,\\
&s_2=\theta(s,\mu),\\
&s_3=\eta(s,\mu),\\
&s_4=\zeta(s,\mu).
\endaligned\right.
\mytag{6.27}
$$
\mytheorem{6.1} The functions \mythetag{6.27}, where $\theta(s,\mu)$, $\eta(s,\mu)$, 
and $\zeta(s,\mu)$ are given by the formulas \mythetag{6.24}, \mythetag{6.25} and 
\mythetag{6.26} respectively, provide a two-parametric solution of the slanted
cuboid equation \mythetag{3.9}. 
\endproclaim
    The proof of this theorem is pure computations. These computations are performed 
by means of the following code:
\medskip
\parshape 1 20pt 340pt
\noindent
\darkred{{\tt restart:\newline
u1:=(1/s1-s1)/2: u2:=(1/s2-s2)/2:\newline
u3:=(1/s3-s3)/2: u4:=(1/s4-s4)/2:\newline
u1:=normal(u1):\ \ \ u2:=normal(u2):\newline
u3:=normal(u3):\ \ \ u4:=normal(u4):\newline
\newline
Eq\_5:=2*u1\^{}2+2*u2\^{}2-u3\^{}2-u4\^{}2:\newline
Eq\_5:=numer(normal(Eq\_5));\newline
\newline
theta:=(1-s\^{}2)*((1-mu\^{}2)\^{}2-4*mu\^{}2)/4/mu/s/(1-mu\^{}2);\newline
eta:=4*mu*s*(1-mu\^{}2)/(1-s\^{}2)/(1+mu\^{}2)/(1-mu\^{}2+2*mu);\newline
zeta:=(1-s\^{}2)*(1+mu\^{}2)*(1-mu\^{}2-2*mu)/4/mu/s/(1-mu\^{}2);\newline
\newline
Eq\_5:=subs(s1=s,s2=theta,s3=eta,s4=zeta,Eq\_5):\newline
Eq\_5:=numer(normal(Eq\_5));
}}
\medskip
     Due to \mythetag{3.7} the parameter $s$ in \mythetag{6.27} is restricted 
by the inequalities $0<s<1$. The parameter $m$ in \mythetag{4.3} is restricted
by the inequalities $0<m<1$. But due to the factor $1-\mu^2-2\,\mu$ in
\mythetag{6.26} and $s_4>0$ in \mythetag{3.7} we have an auxiliary restriction:
$$
\hskip -2em
1-\mu^2-2\,\mu>0.
\mytag{6.28}
$$
Resolving \mythetag{6.28} with respect to $\mu$, we get $\mu<\sqrt{2}-1$. Therefore
the functions $\theta(s,\mu)$, $\eta(s,\mu)$, and $\zeta(s,\mu)$ are well-defined for
$$
\xalignat 2
&\hskip -2em
0<s<1,
&&0<\mu<\sqrt{2}-1.
\mytag{6.29}
\endxalignat 
$$
But due to the inequalities \mythetag{3.7} and \mythetag{3.10} the actual domain 
$D$ of the mapping \mythetag{6.27} could be even smaller than \mythetag{6.29}.\par
     The image of the domain $D$ under the mapping \mythetag{6.24} is a two-dimensional
real algebraic subvariety within $\Gamma_{3\sssize ++}$. Above we have denoted it
through $\Gamma_2^1$. Note that the slanted cuboid equation \mythetag{3.9} and 
the slanted cuboid inequalities \mythetag{3.7} and \mythetag{3.10} admit the following
two discrete symmetry transformations:
$$
\xalignat 2
&\hskip -2em
s_1\longleftrightarrow s_2,
&&s_3\longleftrightarrow s_4
\mytag{6.30}
\endxalignat 
$$
Applying \mythetag{6.30} to \mythetag{6.27}, we derive three more mappings:
$$
\xalignat 3
&\hskip -2em
\left\{
\aligned
&s_1=\theta(s,\mu),\\
&s_2=s,\\
&s_3=\eta(s,\mu),\\
&s_4=\zeta(s,\mu),
\endaligned\right.
&&\left\{
\aligned
&s_1=s,\\
&s_2=\theta(s,\mu),\\
&s_3=\zeta(s,\mu),\\
&s_4=\eta(s,\mu),
\endaligned\right.
&&\left\{
\aligned
&s_1=\theta(s,\mu),\\
&s_2=s,\\
&s_3=\zeta(s,\mu),\\
&s_4=\eta(s,\mu).
\endaligned\right.
\quad
\mytag{6.31}
\endxalignat 
$$
The images of the domain $D$ under the mappings \mythetag{6.31} constitute
three more two-dimensional real algebraic subvarieties within 
$\Gamma_{3\sssize ++}$. We denote them $\Gamma_2^2$, $\Gamma_2^3$, and
$\Gamma_2^4$ respectively.\par
     In Walter Wyss's paper \mycite{4} we find two examples of rational slanted
cuboids. The first example on page 14 is produced by choosing 
$$
\xalignat 2
&\hskip -2em
s=\frac{1}{2},
&&\mu=\frac{1}{3}.
\mytag{6.32}
\endxalignat
$$
Substituting \mythetag{6.32} into \mythetag{6.27} and taking into account 
\mythetag{6.24}, \mythetag{6.25}, \mythetag{6.26}, we get
$$
\xalignat 4
&\hskip -2em
s_1=\frac{1}{2},
&&s_2=\frac{7}{16},
&&s_3=\frac{16}{35},
&&s_4=\frac{5}{16}.
\quad
\mytag{6.33}
\endxalignat
$$
The values \mythetag{6.33} are produced by the following code:
\medskip
\parshape 1 20pt 340pt
\noindent
\darkred{{\tt s1=subs(s=1/2,mu=1/3,s);\newline
s2=subs(s=1/2,mu=1/3,theta);\newline
s3=subs(s=1/2,mu=1/3,eta);\newline
s4=subs(s=1/2,mu=1/3,zeta);
}}
\medskip
\noindent
They do coincide with Walter Wiss's data on page 14. The second example on page
15 of \mycite{4} is produced by choosing
$$
\xalignat 2
&\hskip -2em
s=\frac{12}{25},
&&\mu=\frac{1}{3}.
\mytag{6.34}
\endxalignat
$$
Substituting \mythetag{6.34} into \mythetag{6.27} and taking into account 
\mythetag{6.24}, \mythetag{6.25}, \mythetag{6.26}, we get
$$
\xalignat 4
&\hskip -2em
s_1=\frac{12}{25},
&&s_2=\frac{3367}{7200},
&&s_3=\frac{1440}{3367},
&&s_4=\frac{481}{1440}.
\quad
\mytag{6.35}
\endxalignat
$$
The values \mythetag{6.35} again coincide with Walter Wiss's data on page 15
of his paper.
\head
7. Further verifications.
\endhead
    In sections 6 and 7 of his paper \mycite{4} Walter Wyss changes some notations.
Nevertheless we can continue verifying his formulas relying on Theorem~\mythetheorem{3.1}
and referring them to the basic equation \mythetag{3.9} of the slanted cuboids. The basic
equation \mythetag{3.9} is programmed by means of the following code:
\medskip
\parshape 1 20pt 340pt
\noindent
\darkred{{\tt restart:\newline
v1:=(s1+1/s1)/2: v2:=(s2+1/s2)/2:\newline
v3:=(s3+1/s3)/2: v4:=(s4+1/s4)/2:
}}
\medskip
\parshape 1 20pt 340pt
\noindent
\darkred{{\tt u1:=(1/s1-s1)/2: u2:=(1/s2-s2)/2:\newline
u3:=(1/s3-s3)/2: u4:=(1/s4-s4)/2:\newline 
v1:=normal(v1):  v2:=normal(v2):\newline
v3:=normal(v3):  v4:=normal(v4):\newline
u1:=normal(u1):  u2:=normal(u2):\newline
u3:=normal(u3):  u4:=normal(u4):\newline
\newline
Eq\_5:=2*u1\^{}2+2*u2\^{}2-u3\^{}2-u4\^{}2:\newline
Eq\_5:=numer(normal(Eq\_5));
}}
\medskip
In the beginning of section 7 of his paper \mycite{4} on page 18 Walter Wyss
writes 5 equations which are not numbered. Two of them coincide with the
equations (16) and (17) on page 3. Other three equations coincide with the
equations (9), (10), (11) on page 2. The equation (8) from page 2 is not 
written on page 18. Like in the case of the equations (16), (18), and (19)
on page 3, we have a subset of the slanted cuboid equations \mythetag{3.4}
and \mythetag{3.5}. They are fulfilled once the basic equation \mythetag{3.9}
is fulfilled and the formulas \mythetag{3.6} are taken for expressing
$u_1$, $u_2$, $u_3$, $u_4$ and $v_1$, $v_2$, $v_3$, $v_4$ through the generators
$s_1$, $s_2$, $s_3$, $s_4$. The numbers \mythetag{5.6} and \mythetag{5.7} are
expressed on page 16 of Walter Wyss's paper \mycite{4} in a functional form as
values of some function $m(x)$ (see formulas (95)). We define this
function as
\medskip
\parshape 1 20pt 340pt
\noindent
\darkred{{\tt m:=proc(x) option remember: end proc:
}}
\medskip
\noindent
Its values $m(\alpha)$ and $m(\alpha_1)$ are coded as follows:
\medskip
\parshape 1 20pt 340pt
\noindent
\darkred{{\tt m(alpha):=normal((2*u2+u3-u4)/(2*u1+u3+u4)):\newline
m(alpha1):=normal((2*v2+v3-v4)/(2*u1+v3+v4)):
}}
\medskip
\noindent
The formulas (96) from Walter Wyss's paper \mycite{4} are coded similarly:
\medskip
\parshape 1 20pt 340pt
\noindent
\darkred{{\tt m(beta):=normal((2*u2-u3+u4)/(2*u1+u3+u4)):\newline
m(beta1):=normal((2*v2-v3+v4)/(2*u1+v3+v4)):
}}
\medskip
\noindent
After the formulas (96) on page 16 we see some formulas which are not numbered 
Some of them coincide with the not numbered formulas on page 4. They lead
to \mythetag{5.12} and \mythetag{5.13}, where $\sigma_1$ is defined by means
of the formula
$$
\hskip -2em
\sigma_1:=\frac{\alpha+\alpha_1}{2}
\mytag{7.1}
$$
Trigonometric functions of the angle $\alpha$ and its multiples are coded as
follows:
\medskip
\parshape 1 20pt 340pt
\noindent
\darkred{{\tt unprotect(sin,cos,tan,cot):\newline
sin:=subsop(3=NULL,eval(sin)):\newline
cos:=subsop(3=NULL,eval(cos)):\newline
tan:=subsop(3=NULL,eval(tan)):\newline
cot:=subsop(3=NULL,eval(cot)):\newline
\newline
tan(alpha/2):=m(alpha):\newline
cot(alpha/2):=1/m(alpha):\newline
sin(alpha/2):=m(alpha)*cos(alpha/2):\newline
sin(alpha):=normal(2*m(alpha)/(1+m(alpha)\^{}2)):\newline
cos(alpha):=normal((1-m(alpha)\^{}2)/(1+m(alpha)\^{}2)):\newline
}}
\medskip
\noindent
Trigonometric functions of the angle $\alpha_1$ and its multiples are coded 
similarly:
\medskip
\parshape 1 20pt 340pt
\noindent
\darkred{{\tt tan(alpha1/2):=m(alpha1):\newline
cot(alpha1/2):=1/m(alpha1):\newline
sin(alpha1/2):=m(alpha1)*cos(alpha1/2):\newline
sin(alpha1):=normal(2*m(alpha1)/(1+m(alpha1)\^{}2)):\newline
cos(alpha1):=normal((1-m(alpha1)\^{}2)/(1+m(alpha1)\^{}2)):
}}
\medskip
\noindent
The formulas \mythetag{5.13} and \mythetag{7.1} are coded as follows:
\medskip
\parshape 1 20pt 340pt
\noindent
\darkred{{\tt psi:=Pi/4-alpha/2-alpha1/2:\newline
sigma1:=alpha/2+alpha1/2:
}}
\medskip
\noindent
The functions $\omega_{+}(x)$ and $\omega_{-}(x)$ are defined according to
\mythetag{4.10}:
\medskip
\parshape 1 20pt 340pt
\noindent
\darkred{{\tt omega\_plus:=proc(x) sin(x)+cos(x) end proc:\newline
omega\_minus:=proc(x) cos(x)-sin(x) end proc:
}}
\medskip
\noindent
Now we are able to verify the formulas \mythetag{5.12} which are repeated on
page 16 of Walter Wyss's paper \mycite{4}. This is done by means of the following
code:
\medskip
\parshape 1 20pt 340pt
\noindent
\darkred{{\tt Eq\_5\_12\_1:=omega\_plus(sigma1)-sqrt(2)*cos(psi):\newline
Eq\_5\_12\_1:=expand(Eq\_5\_12\_1):\newline
Eq\_5\_12\_2:=omega\_minus(sigma1)-sqrt(2)*sin(psi):\newline
Eq\_5\_12\_2:=expand(Eq\_5\_12\_2):\newline
\newline
Eq\_5\_12\_1,Eq\_5\_12\_2;
}}
\medskip
\noindent
The following formula on page 16 of the paper \mycite{4} is immediate from 
\mythetag{5.12}:
$$
\hskip -2em
\lambda=\tan\psi=\frac{\omega_{-}(\sigma_1)}{\omega_{+}(\sigma_1)}.
\mytag{7.2}
$$
The equation \mythetag{7.2} can be verified directly by means of the following code:
\medskip
\parshape 1 20pt 340pt
\noindent
\darkred{{\tt Eq\_7\_2:=tan(psi)-omega\_minus(sigma1)/omega\_plus(sigma1):\newline
Eq\_7\_2:=expand(Eq\_7\_2):\newline
Eq\_7\_2:=numer(normal(Eq\_7\_2)):\newline
Eq\_7\_2:=rem(Eq\_7\_2,Eq\_5,s1);
}}
\medskip
     The formula (97) in paper \mycite{4} is just a notation. It is coded as follows:
\medskip
\parshape 1 20pt 340pt
\noindent
\darkred{{\tt k:=(m(alpha)+m(alpha1))/(1-m(alpha)*m(alpha1)):\newline
k:=normal(k):
}}
\medskip
\noindent
The formulas (98) on page 17 of the paper \mycite{4} are different. They should be 
verified. We verify them by means of the following code:
\medskip
\parshape 1 20pt 340pt
\noindent
\darkred{{\tt Eq\_98\_1:=sin(2*sigma1)-2*k/(1+k\^{}2):\newline
Eq\_98\_1:=expand(Eq\_98\_1):\newline
Eq\_98\_1:=numer(normal(Eq\_98\_1)):\newline
Eq\_98\_1:=rem(Eq\_98\_1,Eq\_5,s1):\newline
\newline
Eq\_98\_2:=cos(2*sigma1)-(1-k\^{}2)/(1+k\^{}2):\newline
Eq\_98\_2:=expand(Eq\_98\_2):\newline
Eq\_98\_2:=numer(normal(Eq\_98\_2)):
}}
\medskip
\parshape 1 20pt 340pt
\noindent
\darkred{{\tt Eq\_98\_2:=rem(Eq\_98\_2,Eq\_5,s1):\newline
\newline
Eq\_98\_1,Eq\_98\_2;
}}
\medskip
\noindent
The formula (99) is immediate from the formulas \mythetag{6.9} and \mythetag{6.10},
which are the identities with respect to the argument $\alpha$.\par
     The formula (100) in \mycite{4} follows from (98) and (99). But, nevertheless,
we verify this formula directly by means of the following code:
\medskip
\parshape 1 20pt 340pt
\noindent
\darkred{{\tt Eq\_100:=tan(psi)-(1-k)/(1+k):\newline
Eq\_100:=expand(Eq\_100):\newline
Eq\_100:=numer(normal(Eq\_100)):\newline
Eq\_100:=rem(Eq\_100,Eq\_5,s1);
}}
\medskip
\noindent
The formulas (101) and (102) in \mycite{4} are written by analogy to the formulas
(97) and (100). These formulas are coded as follows:
\medskip
\parshape 1 20pt 340pt
\noindent
\darkred{{\tt bar\_k:=(m(beta)+m(beta1))/(1-m(beta)*m(beta1)):\newline
bar\_lambda:=(1-bar\_k)/(1+bar\_k):
}}
\medskip
\noindent
On page 17 of Walter Wyss's paper \mycite{4} we see the phrase: ``Therefore the 
parameters $u_1$, $\beta$, $\bar\lambda$ also satisfy the general equations, 
however with the interchange of $s_3$ with $s_4$''. We do not know which general
equations does he mean. But we suspect that he means the equations (82), (83), (84)
on page 11 of his paper. Their analogs for the variables $u_1$, $\beta$, $\bar\lambda$ 
look like
$$
\gather
\hskip -2em
\omega_{-}(\beta)-\bar\lambda\,\omega_{+}(\beta)
=s_2\,s_4+\bar\lambda\,s_2\,s_3,
\mytag{7.3}\\
\hskip -2em
Q\,(\omega_{+}(\beta)+\bar\lambda\,\omega_{-}(\beta))=s_2\,s_4-\bar\lambda\,s_2\,s_3,
\mytag{7.4}\\
\hskip -2em
2\,u_1\,(\omega_{-}(\beta)-\bar\lambda\,\omega_{+}(\beta))+s_3-\bar\lambda\,s_4
=s_2\,(\omega_{+}(\beta)+\bar\lambda\,\omega_{-}(\beta)).
\mytag{7.5}
\endgather
$$
Here $Q$ is given by the formula \mythetag{5.10} and $\bar\lambda$ is given
by the formula (102) in \mycite{4}. The formulas \mythetag{7.3}, 
\mythetag{7.4}, \mythetag{7.5} are verified as follows:
\medskip
\parshape 1 20pt 340pt
\noindent
\darkred{{\tt tan(beta/2):=m(beta):\newline
cot(beta/2):=1/m(beta):\newline
sin(beta/2):=m(beta)*cos(beta/2):\newline
sin(beta):=normal(2*m(beta)/(1+m(beta)\^{}2)):\newline
cos(beta):=normal((1-m(beta)\^{}2)/(1+m(beta)\^{}2)):\newline
\newline
tan(beta1/2):=m(beta1):\newline
cot(beta1/2):=1/m(beta1):\newline
sin(beta1/2):=m(beta1)*cos(beta1/2):\newline
sin(beta1):=normal(2*m(beta1)/(1+m(beta1)\^{}2)):\newline
cos(beta1):=normal((1-m(beta1)\^{}2)/(1+m(beta1)\^{}2)):\newline
\newline
Eq\_7\_3:=-omega\_minus(beta)+bar\_lambda*omega\_plus(beta)\newline
\phantom{aaaaaaaa}+s2*s4+bar\_lambda*s2*s3:\newline
Eq\_7\_3:=expand(Eq\_7\_3):\newline
Eq\_7\_3:=numer(normal(Eq\_7\_3)):\newline
Eq\_7\_3:=rem(Eq\_7\_3,Eq\_5,s1):\newline
}}
\medskip
\parshape 1 20pt 340pt
\noindent
\darkred{{\tt Q:=s3*s4:\newline
Eq\_7\_4:=-Q*(omega\_plus(beta)+bar\_lambda*omega\_minus(beta))\newline
\phantom{aaaaaaaa}+s2*s4-bar\_lambda*s2*s3:\newline
Eq\_7\_4:=expand(Eq\_7\_4):\newline
Eq\_7\_4:=numer(normal(Eq\_7\_4)):\newline
Eq\_7\_4:=rem(Eq\_7\_4,Eq\_5,s1):\newline
\newline
Eq\_7\_5:=-2*u1*(omega\_minus(beta)-bar\_lambda*omega\_plus(beta))\newline
\phantom{aaaaaaaa}-s3+bar\_lambda*s4+s2*(omega\_plus(beta)\newline
\phantom{aaaaaaaa}+bar\_lambda*omega\_minus(beta)):\newline
Eq\_7\_5:=expand(Eq\_7\_5):\newline
Eq\_7\_5:=numer(normal(Eq\_7\_5)):\newline
Eq\_7\_5:=rem(Eq\_7\_5,Eq\_5,s1):\newline
\newline
Eq\_7\_3,Eq\_7\_4,Eq\_7\_5;
}}
\medskip
     In the beginning of section 7 of his paper \mycite{4} Walter Wyss
introduces some new notations replacing previous ones (see (103), (104)
(105), (106) on page 18): 
$$
\xalignat 2
&\hskip -2em
\sigma=\frac{\alpha+\beta}{2},
&&\delta=\frac{\alpha-\beta}{2}.
\mytag{7.6}\\
&\hskip -2em
\sigma_1=\frac{\alpha_1+\beta_1}{2},
&&\delta_1=\frac{\alpha_1-\beta_1}{2}.
\mytag{7.7}\\
\endxalignat
$$
The notations \mythetag{7.7} replace the notations introduced on page 4
of the paper \mycite{4}. The second notation \mythetag{7.6} replaces the 
notation used on page 9 of the paper \mycite{4}. The notations \mythetag{7.6} 
do coincide with the notations \mythetag{4.9}. As a whole, the notations
\mythetag{7.6} and \mythetag{7.7} are coded as follows:
\medskip
\parshape 1 20pt 340pt
\noindent
\darkred{{\tt sigma:=alpha/2+beta/2:\newline
delta:=alpha/2-beta/2:\newline
sigma1:=alpha1/2+beta1/2:\newline
delta1:=alpha1/2-beta1/2:
}}
\medskip
    The next are the formulas (107), (108). The formulas (107) are verified
as follows:
\medskip
\parshape 1 20pt 340pt
\noindent
\darkred{{\tt Eq\_107\_1:=u2-u1*tan(sigma):\newline
Eq\_107\_1:=expand(Eq\_107\_1):\newline
Eq\_107\_1:=numer(normal(Eq\_107\_1)):\newline
Eq\_107\_1:=rem(Eq\_107\_1,Eq\_5,s1):\newline
\newline
Eq\_107\_2:=u3-u1*omega\_plus(delta)/cos(sigma):\newline
Eq\_107\_2:=expand(Eq\_107\_2):\newline
Eq\_107\_2:=numer(normal(Eq\_107\_2)):\newline
Eq\_107\_2:=rem(Eq\_107\_2,Eq\_5,s1):\newline
\newline
Eq\_107\_3:=u4-u1*omega\_minus(delta)/cos(sigma):\newline
Eq\_107\_3:=expand(Eq\_107\_3):\newline
Eq\_107\_3:=numer(normal(Eq\_107\_3)):\newline
Eq\_107\_3:=rem(Eq\_107\_3,Eq\_5,s1):\newline
\newline
Eq\_107\_1,Eq\_107\_2,Eq\_107\_3;
}}
\medskip
\noindent
The formulas (108) in \mycite{4} are verified similarly:
\medskip
\parshape 1 20pt 340pt
\noindent
\darkred{{\tt Eq\_108\_1:=v2-u1*tan(sigma1):\newline
Eq\_108\_1:=expand(Eq\_108\_1):\newline
Eq\_108\_1:=numer(normal(Eq\_108\_1)):\newline
Eq\_108\_1:=rem(Eq\_108\_1,Eq\_5,s1):\newline
\newline
Eq\_108\_2:=v3-u1*omega\_plus(delta1)/cos(sigma1):\newline
Eq\_108\_2:=expand(Eq\_108\_2):\newline
Eq\_108\_2:=numer(normal(Eq\_108\_2)):\newline
Eq\_108\_2:=rem(Eq\_108\_2,Eq\_5,s1):\newline
\newline
Eq\_108\_3:=v4-u1*omega\_minus(delta1)/cos(sigma1):\newline
Eq\_108\_3:=expand(Eq\_108\_3):\newline
Eq\_108\_3:=numer(normal(Eq\_108\_3)):\newline
Eq\_108\_3:=rem(Eq\_108\_3,Eq\_5,s1):\newline
\newline
Eq\_108\_1,Eq\_108\_2,Eq\_108\_3;
}}
\medskip
\noindent
The formulas (109) and (110) are converse to (107) and (108). Nevertheless, we
verify them directly with the use of the following code:
\medskip
\parshape 1 20pt 340pt
\noindent
\darkred{{\tt Eq\_109\_1:=tan(sigma)-u2/u1:\newline
Eq\_109\_1:=expand(Eq\_109\_1):\newline
Eq\_109\_1:=numer(normal(Eq\_109\_1)):\newline
Eq\_109\_1:=rem(Eq\_109\_1,Eq\_5,s1):\newline
\newline
Eq\_109\_2:=tan(delta)-(u3-u4)/(u3+u4):\newline
Eq\_109\_2:=expand(Eq\_109\_2):\newline
Eq\_109\_2:=numer(normal(Eq\_109\_2)):\newline
Eq\_109\_2:=rem(Eq\_109\_2,Eq\_5,s1):\newline
\newline
Eq\_110\_1:=tan(sigma1)-v2/u1:\newline
Eq\_110\_1:=expand(Eq\_110\_1):\newline
Eq\_110\_1:=numer(normal(Eq\_110\_1)):\newline
Eq\_110\_1:=rem(Eq\_110\_1,Eq\_5,s1):\newline
\newline
Eq\_110\_2:=tan(delta1)-(v3-v4)/(v3+v4):\newline
Eq\_110\_2:=expand(Eq\_110\_2):\newline
Eq\_110\_2:=numer(normal(Eq\_110\_2)):\newline
Eq\_110\_2:=rem(Eq\_110\_2,Eq\_5,s1):\newline
\newline
Eq\_109\_1,Eq\_109\_2,Eq\_110\_1,Eq\_110\_2;
}}
\medskip
    The next are the formulas (111), (112), and (113). In the case of the first
formula (111) we use the following code in order to verify it:
\medskip
\parshape 1 20pt 340pt
\noindent
\darkred{{\tt Eq\_111\_1:=u1\^{}2*(tan(sigma1)\^{}2-tan(sigma)\^{}2)-1:\newline
Eq\_111\_1:=expand(Eq\_111\_1):\newline
Eq\_111\_1:=numer(normal(Eq\_111\_1)):\newline
Eq\_111\_1:=rem(Eq\_111\_1,Eq\_5,s1);
}}
\medskip
\noindent
The second formula (111) is more complicated for computer handling:
\medskip
\parshape 1 20pt 340pt
\noindent
\darkred{{\tt psi\_phi\_simplify:=proc(A) local AA: global m:\newline
\phantom{aa}AA:=subs(cos(alpha/2)\^{}2=1/(1+m(alpha)\^{}2),A):\newline
\phantom{aa}AA:=subs(cos(alpha1/2)\^{}2=1/(1+m(alpha1)\^{}2),AA):\newline
\phantom{aa}AA:=subs(cos(beta/2)\^{}2=1/(1+m(beta)\^{}2),AA):\newline
\phantom{aa}AA:=subs(cos(beta1/2)\^{}2=1/(1+m(beta1)\^{}2),AA):\newline
\phantom{aa}return AA:\newline
end proc:\newline
\newline
uuu:=psi\_phi\_simplify(expand(cos(sigma)\^{}2)):\newline
uuu1:=psi\_phi\_simplify(expand(cos(sigma1)\^{}2)):\newline
Eq\_111\_2:=u1\^{}2*(1/uuu1-1/uuu)-1:\newline
Eq\_111\_2:=expand(Eq\_111\_2):\newline
Eq\_111\_2:=numer(normal(Eq\_111\_2)):\newline
Eq\_111\_2:=rem(Eq\_111\_2,Eq\_5,s1);
}}
\medskip
\noindent
Now let's proceed to the formulas (112). For them we use the following code:
\medskip
\parshape 1 20pt 340pt
\noindent
\darkred{{\tt vvv:=psi\_phi\_simplify(expand(omega\_plus(delta)\^{}2)):\newline
vvv1:=psi\_phi\_simplify(expand(omega\_plus(delta1)\^{}2)):\newline
Eq\_112\_1:=u1\^{}2*(vvv1/uuu1-vvv/uuu)-1:\newline
Eq\_112\_1:=numer(normal(Eq\_112\_1)):\newline
Eq\_112\_1:=rem(Eq\_112\_1,Eq\_5,s1);\newline
\newline
vvv:=expand(1+sin(2*delta)):\newline
vvv1:=expand(1+sin(2*delta1)):\newline
Eq\_112\_2:=u1\^{}2*(vvv1/uuu1-vvv/uuu)-1:\newline
Eq\_112\_2:=numer(normal(Eq\_112\_2)):\newline
Eq\_112\_2:=rem(Eq\_112\_2,Eq\_5,s1);
}}
\medskip
\noindent
The formulas (113) are similar to (112). For them we use the following code:
\medskip
\parshape 1 20pt 340pt
\noindent
\darkred{{\tt vvv:=psi\_phi\_simplify(expand(omega\_minus(delta)\^{}2)):\newline
vvv1:=psi\_phi\_simplify(expand(omega\_minus(delta1)\^{}2)):\newline
Eq\_113\_1:=u1\^{}2*(vvv1/uuu1-vvv/uuu)-1:\newline
Eq\_113\_1:=numer(normal(Eq\_113\_1)):\newline
Eq\_113\_1:=rem(Eq\_113\_1,Eq\_5,s1);\newline
\newline
vvv:=expand(1-sin(2*delta)):\newline
vvv1:=expand(1-sin(2*delta1)):\newline
Eq\_113\_2:=u1\^{}2*(vvv1/uuu1-vvv/uuu)-1:\newline
Eq\_113\_2:=numer(normal(Eq\_113\_2)):\newline
Eq\_113\_2:=rem(Eq\_113\_2,Eq\_5,s1);
}}
\medskip
     The next are the formulas (114) and (115). They are written as follows:
$$
\gather
\hskip -2em
\tan^2\sigma_1-\tan^2\sigma=\tan^2\psi,
\mytag{7.8}\\
\vspace{1ex}
\hskip -2em
\frac{sin(2\,\delta)}{\cos^2\sigma}=\frac{sin(2\,\delta_1)}{\cos^2\sigma_1}.
\mytag{7.9}
\endgather
$$
{\bf Attention!} At the bottom of page 18 in his paper \mycite{4} Walter Wyss
writes: ``We rename $\psi_1=\psi$''. This means that in \mythetag{7.8} $\psi$ 
does not coincide with \mythetag{5.13}. It coincides with
$\psi_1$ in \mythetag{3.8}. From \mythetag{3.8} we derive
$$
\hskip -2em
\tan\psi=\tan\psi_1=\frac{1}{u_1}. 
\mytag{7.10}
$$
We use \mythetag{7.10} in writing code for verifying the formula \mythetag{7.8}:
\medskip
\parshape 1 20pt 340pt
\noindent
\darkred{{\tt Eq\_7\_8:=tan(sigma1)\^{}2-tan(sigma)\^{}2-1/u1\^{}2:\newline
Eq\_7\_8:=expand(Eq\_7\_8):\newline
Eq\_7\_8:=numer(normal(Eq\_7\_8)):\newline
Eq\_7\_8:=rem(Eq\_7\_8,Eq\_5,s1);
}}
\medskip
\noindent
In the case of the formula \mythetag{7.9} we use the following code:
\medskip
\parshape 1 20pt 340pt
\noindent
\darkred{{\tt vvv:=expand(sin(2*delta)):\newline
vvv1:=expand(sin(2*delta1)):\newline
Eq\_7\_9:=vvv1/uuu1-vvv/uuu:\newline
Eq\_7\_9:=numer(normal(Eq\_7\_9)):\newline
Eq\_7\_9:=rem(Eq\_7\_9,Eq\_5,s1);
}}
\medskip
    The formulas (116), (117), (118), and (119) in Walter Wyss's paper \mycite{4}
are just notations. They are coded as follows:
\medskip
\parshape 1 20pt 340pt
\noindent
\darkred{{\tt M:=normal(expand(tan(sigma))):\newline
M1:=normal(expand(tan(sigma1))):\newline
N:=normal(expand(tan(delta))):\newline
N1:=normal(expand(tan(delta1))):
}}
\medskip
\noindent
The formulas (120) and (121) are trivial. Nevertheless, we can verify them:
\medskip
\parshape 1 20pt 340pt
\noindent
\darkred{{\tt Eq\_120:=sin(alpha)/cos(alpha)-(M+N)/(1-M*N):\newline
Eq\_120:=numer(normal(Eq\_120)):\newline
Eq\_120:=rem(Eq\_120,Eq\_5,s1):\newline
\newline
Eq\_121:=sin(alpha1)/cos(alpha1)-(M1+N1)/(1-M1*N1):\newline
Eq\_121:=numer(normal(Eq\_121)):\newline
Eq\_121:=rem(Eq\_121,Eq\_5,s1):\newline
\newline
Eq\_120, Eq\_121;
}}
\medskip
\noindent
The formulas (122) and (123) follow from (120) and (121). But we
verify them too:
\medskip
\parshape 1 20pt 340pt
\noindent
\darkred{{\tt Eq\_122:=N-(sin(alpha)/cos(alpha)-M)\newline
\phantom{aaaaaaa}/(1+M*sin(alpha)/cos(alpha)):\newline
Eq\_122:=numer(normal(Eq\_122)):\newline
Eq\_122:=rem(Eq\_122,Eq\_5,s1):\newline
Eq\_123:=N1-(sin(alpha1)/cos(alpha1)-M1)\newline
\phantom{aaaaaaa}/(1+M1*sin(alpha1)/cos(alpha1)):\newline
Eq\_123:=numer(normal(Eq\_123)):\newline
Eq\_123:=rem(Eq\_123,Eq\_5,s1):\newline
\newline
Eq\_122, Eq\_123;
}}
\medskip
\noindent
Taking into account \mythetag{7.10}, we can verify the formula (124) as follows:
\medskip
\parshape 1 20pt 340pt
\noindent
\darkred{{\tt Eq\_124:=M1\^{}2-M\^{}2-1/u1\^{}2:\newline
Eq\_124:=numer(normal(Eq\_124)):\newline
Eq\_124:=rem(Eq\_124,Eq\_5,s1);
}}
\medskip
\noindent
Then we verify the formula (125) by means of the following code:
\medskip
\parshape 1 20pt 340pt
\noindent
\darkred{{\tt Eq\_125:=2*N*(1+M\^{}2)/(1+N\^{}2)-2*N1*(1+M1\^{}2)/(1+N1\^{}2):\newline
Eq\_125:=numer(normal(Eq\_125)):\newline
Eq\_125:=rem(Eq\_125,Eq\_5,s1);
}}
\medskip
\noindent
The formula (126) is verified similarly. We use the following code for it:
\medskip
\parshape 1 20pt 340pt
\noindent
\darkred{{\tt vvv:=normal(expand(sin(2*alpha))):\newline
vvv1:=normal(expand(sin(2*alpha1))):\newline
uuu:=normal(expand(cos(2*alpha))):\newline
uuu1:=normal(expand(cos(2*alpha1))):\newline
Eq\_126:=(1-M\^{}2)*vvv-2*M*uuu-(1-M1\^{}2)*vvv1+2*M1*uuu1:\newline
Eq\_126:=numer(normal(Eq\_126)):\newline
Eq\_126:=rem(Eq\_126,Eq\_5,s1);
}}
\medskip
     The formulas (127), (128), (129) are elementary. They follow from the 
notations \mythetag{7.6} and \mythetag{7.7} complemented with the notation
(116). So we proceed to the formulas (130) and (131). The formula (130)
is verified with the use of the code
\medskip
\parshape 1 20pt 340pt
\noindent
\darkred{{\tt vvv:=factor(expand(omega\_plus(delta))):\newline
uuu:=factor(expand(cos(sigma))):\newline
www:=omega\_plus(alpha)+omega\_minus(beta):\newline
Eq\_130:=normal(vvv/uuu)-1/2*(1+M\^{}2)*www:\newline
Eq\_130:=numer(normal(Eq\_130)):\newline
Eq\_130:=rem(Eq\_130,Eq\_5,s1);
}}
\medskip
\noindent
The formula (131) is similar. In this case we use the following code:
\medskip
\parshape 1 20pt 340pt
\noindent
\darkred{{\tt vvv:=factor(expand(omega\_minus(delta))):\newline
uuu:=factor(expand(cos(sigma))):\newline
www:=omega\_minus(alpha)+omega\_plus(beta):\newline
Eq\_131:=normal(vvv/uuu)-1/2*(1+M\^{}2)*www:\newline
Eq\_131:=numer(normal(Eq\_131)):\newline
Eq\_131:=rem(Eq\_131,Eq\_5,s1);
}}
\medskip
    The formulas (132) and (134) are elementary. Therefore we proceed to
the formulas (133) and (135). They are verified with the following code:
\medskip
\parshape 1 20pt 340pt
\noindent
\darkred{{\tt Eq\_133:=cos(beta)-(1-M\^{}2)/(1+M\^{}2)*cos(alpha)\newline
\phantom{aaaaaaa}-2*M/(1+M\^{}2)*sin(alpha):\newline
Eq\_133:=numer(normal(Eq\_133)):\newline
Eq\_133:=rem(Eq\_133,Eq\_5,s1);\newline
\newline
Eq\_135:=sin(beta)-2*M/(1+M\^{}2)*cos(alpha)\newline
\phantom{aaaaaaa}+(1-M\^{}2)/(1+M\^{}2)*sin(alpha):\newline
Eq\_135:=numer(normal(Eq\_135)):\newline
Eq\_135:=rem(Eq\_135,Eq\_5,s1);
}}
\medskip
    The next are the formulas (136) and (137). The first of these two formulas 
is verified by means of the following code:
\medskip
\parshape 1 20pt 340pt
\noindent
\darkred{{\tt Eq\_136:=omega\_plus(beta)-(1-M\^{}2)/(1+M\^{}2)*omega\_minus(alpha)\newline
\phantom{aaaaaaa}-2*M/(1+M\^{}2)*omega\_plus(alpha):\newline
Eq\_136:=numer(normal(Eq\_136)):\newline
Eq\_136:=rem(Eq\_136,Eq\_5,s1);
}}
\medskip
\noindent
The second one is similar. It is verified with the use of the code
\medskip
\parshape 1 20pt 340pt
\noindent
\darkred{{\tt Eq\_137:=omega\_minus(beta)-(1-M\^{}2)/(1+M\^{}2)*omega\_plus(alpha)\newline
\phantom{aaaaaaa}+2*M/(1+M\^{}2)*omega\_minus(alpha):\newline
Eq\_137:=numer(normal(Eq\_137)):\newline
Eq\_137:=rem(Eq\_137,Eq\_5,s1);
}}
\medskip
\noindent
The formulas (138) and (139) are similar to the previous two formulas (136) and (137). 
They are verified as follows:
\medskip
\parshape 1 20pt 340pt
\noindent
\darkred{{\tt Eq\_138:=omega\_plus(alpha)+omega\_minus(beta)\newline
\phantom{aaaaaaa}-2/(1+M\^{}2)*(omega\_plus(alpha)-M*omega\_minus(alpha)):\newline
Eq\_138:=numer(normal(Eq\_138)):\newline
Eq\_138:=rem(Eq\_138,Eq\_5,s1);\newline
\newline
Eq\_139:=omega\_minus(alpha)+omega\_plus(beta)\newline
\phantom{aaaaaaa}-2/(1+M\^{}2)*(omega\_minus(alpha)+M*omega\_plus(alpha)):\newline
Eq\_139:=numer(normal(Eq\_139)):\newline
Eq\_139:=rem(Eq\_139,Eq\_5,s1);
}}
\medskip
     The next are the formulas (140), (141), (142). They look more simple than
the previous ones. We verify these formulas as follows:
\medskip
\parshape 1 20pt 340pt
\noindent
\darkred{{\tt Eq\_140:=u2-u1*M:\newline
Eq\_140:=numer(normal(Eq\_140)):\newline
Eq\_140:=rem(Eq\_140,Eq\_5,s1);\newline
\newline
Eq\_141:=u3-u1*(omega\_plus(alpha)-M*omega\_minus(alpha)):\newline
Eq\_141:=numer(normal(Eq\_141)):\newline
Eq\_141:=rem(Eq\_141,Eq\_5,s1);\newline
\newline
Eq\_142:=u4-u1*(omega\_minus(alpha)+M*omega\_plus(alpha)):\newline
Eq\_142:=numer(normal(Eq\_142)):\newline
Eq\_142:=rem(Eq\_142,Eq\_5,s1);
}}
\medskip
\noindent
The formulas (143), (144), (145) are similar to the previous formulas
(140), (141), (142). They are verified as follows:
\medskip
\parshape 1 20pt 340pt
\noindent
\darkred{{\tt Eq\_143:=v2-u1*M1:\newline
Eq\_143:=numer(normal(Eq\_143)):\newline
Eq\_143:=rem(Eq\_143,Eq\_5,s1);\newline
\newline
Eq\_144:=v3-u1*(omega\_plus(alpha1)-M1*omega\_minus(alpha1)):\newline
Eq\_144:=numer(normal(Eq\_144)):\newline
Eq\_144:=rem(Eq\_144,Eq\_5,s1);
}}
\medskip
\parshape 1 20pt 340pt
\noindent
\darkred{{\tt Eq\_145:=v4-u1*(omega\_minus(alpha1)+M1*omega\_plus(alpha1)):\newline
Eq\_145:=numer(normal(Eq\_145)):\newline
Eq\_145:=rem(Eq\_145,Eq\_5,s1);
}}
\medskip
\noindent
Exchanging $u_3$ and $u_4$, i\.\,e\. applying the second transformation 
\mythetag{6.30}, Walter Wyss derived two more formulas (146) and (147).
They are verified as follows:
\medskip
\parshape 1 20pt 340pt
\noindent
\darkred{{\tt Eq\_146:=u3-u1*(omega\_minus(beta)+M*omega\_plus(beta)):\newline
Eq\_146:=numer(normal(Eq\_146)):\newline
Eq\_146:=rem(Eq\_146,Eq\_5,s1);\newline
\newline
Eq\_147:=u4-u1*(omega\_plus(beta)-M*omega\_minus(beta)):\newline
Eq\_147:=numer(normal(Eq\_147)):\newline
Eq\_147:=rem(Eq\_147,Eq\_5,s1);
}}
\medskip
     The {\bf cuboid limit} or, being more precise, the {\bf rectangular cuboid limit} 
is the case where the parallelogram ABFE in Fig\.~1.1 turns to a rectangle. In this 
case two its diagonals become equal to each other:
$$
\hskip -2em
|AF|=|EB|. 
\mytag{7.11}
$$
Comparing \mythetag{7.11} with our notations \mythetag{3.3}, we find
$$
\hskip -2em
u_3=u_4.
\mytag{7.12}
$$
Applying \mythetag{3.6} to \mythetag{7.12}, we derive the equation
$$
\hskip -2em
\frac{1-s_3^2}{2\,s_3}=\frac{1-s_4^2}{2\,s_4}.
\mytag{7.13}
$$
The equation \mythetag{7.13} has two solutions
$$
\xalignat 2
&s_3=s_4,
&&s_3=\frac{1}{s_4}.
\endxalignat
$$
But due to the inequalities \mythetag{3.7} only the first solution is suitable
for us:
$$
\hskip -2em
s_3=s_4.
\mytag{7.14}
$$
Thus the {\bf rectangular cuboid limit} is the case where either of the two
equivalent equalities \mythetag{7.12} or \mythetag{7.14} is fulfilled.\par
     On page 21 of his paper \mycite{4} Walter Wyss writes that in the
rectangular cuboid limit the following equalities are fulfilled:
$$
\xalignat 2
&\hskip -2em
N=0,
&&N_1=0.
\mytag{7.15}
\endxalignat
$$
The equalities \mythetag{7.15} are easily verified by means of
the following code:
\medskip
\parshape 1 20pt 340pt
\noindent
\darkred{{\tt subs(s4=s3,N),\newline
subs(s4=s3,N1);
}}
\medskip
\noindent
On page 23 of his paper \mycite{4} Walter Wyss writes that the
cuboid limit is given by 
$$
\xalignat 2
&\hskip -2em
\alpha=\beta,
&&\alpha_1=\beta_1.
\mytag{7.16}
\endxalignat
$$
One can easily verify that the equality \mythetag{7.14} implies both equalities
\mythetag{7.16}. This is done with the use of the following code:
\medskip
\parshape 1 20pt 340pt
\noindent
\darkred{{\tt normal(subs(s4=s3,cos(alpha)-cos(beta))),\newline
normal(subs(s4=s3,sin(alpha)-sin(beta))),\newline
normal(subs(s4=s3,cos(alpha1)-cos(beta1))),\newline
normal(subs(s4=s3,sin(alpha1)-sin(beta1)));
}}
\medskip
\noindent
Apart from \mythetag{7.15} and \mythetag{7.16} there are two more
equalities:
$$
\xalignat 2
&\hskip -2em
M=\tan(\alpha),
&&M_1=\tan(\alpha_1).
\mytag{7.17}
\endxalignat
$$
One can verify that the equality \mythetag{7.14} implies both equalities
\mythetag{7.17}. In this case we do it with the use of the following code:
\medskip
\parshape 1 20pt 340pt
\noindent
\darkred{{\tt subs(s4=s3,normal(M-sin(alpha)/cos(alpha))),\newline
subs(s4=s3,normal(M1-sin(alpha1)/cos(alpha1)));
}}
\medskip
\noindent
Note that the rectangular cuboid limit is not a singular case though 
some Walter Wiss's formulas are not applicable to it. 
\head
8. A special example. 
\endhead
    On page 21 of his paper \mycite{4} Walter Wyss considers a special case 
in the form of an example. This special case is defined by the equality
$$
\hskip -2em
\alpha+\alpha_1=\frac{\pi}{2}.
\mytag{8.1}
$$
On pages 3 and 4 of his paper Walter Wyss writes that $\alpha$, $\alpha_1$, 
and $\alpha_2$ are Heron angles in the first quadrant, i\.\,e\. they obey
the inequalities \mythetag{5.9}. Their generators $m$, $m_1$, and $m_2$ are
given by the formulas (25), (26), and (27) on page 3 of the paper \mycite{4}.
Comparing the formulas (25) and (26) with the formulas (95) on page 16, we
conclude that the angles $\alpha$ and $\alpha_1$ in \mythetag{8.1} are the 
same angles which are used in sections 4 and 5 of Walter Wyss's paper 
\mycite{4}.\par
     Dividing the equality \mythetag{8.1} by $2$, we derive
$$
\hskip -2em
\frac{\alpha+\alpha_1}{2}=\frac{\pi}{4}.
\mytag{8.2}
$$
Substituting \mythetag{8.2} into \mythetag{5.13}, we find that
$$
\hskip -2em
\psi=\frac{\pi}{4}-\frac{\alpha+\alpha_1}{2}=0.
\mytag{8.3}
$$
The equality \mythetag{8.3} implies the equality
$$
\hskip -2em
\tan\psi=0.
\mytag{8.4}
$$
Conversely, applying the inequalities \mythetag{5.9} to \mythetag{5.13}, we
derive the inequality
$$
\hskip -2em
-\frac{\pi}{2}<\psi<\frac{\pi}{2}.
\mytag{8.5}
$$
The tangent function is a monotonic increasing function within the interval 
\mythetag{8.5}. It vanishes exactly once at the point $\psi=0$. This means 
that the equality \mythetag{8.4} implies backward the equality \mythetag{8.3}
and then \mythetag{8.2} and \mythetag{8.1}, i.\,e\. the equalities 
\mythetag{8.1} and \mythetag{8.4} \pagebreak are equivalent.\par
     Now let's recall that the equality \mythetag{8.4} in the form of 
$\lambda=\tan\psi$ and $\lambda=0$ was used by Walter Wyss in order to construct 
a special solution of the slanted cuboid equations (see (81) on page 10 and (89) 
on page 11 of \mycite{4}). Thus the conclusion.
\mytheorem{8.1} The special solution given by the condition \mythetag{6.1} and 
the special example defined by the condition \mythetag{8.1} do coincide. 
\endproclaim
\noindent
{\bf Attention!} Due to renaming variables $\psi_1=\psi$ at the bottom of page 
18 of Walter Wyss's paper \mycite{4} the variable $\psi$ in \mythetag{8.3}, 
\mythetag{8.4}, \mythetag{8.5} does not coincide with this variable on pages
19, 20, 21 and so on.\par
     On page 21 of the paper \mycite{4} we see three formulas, which
are not numbered there:
$$
\gather
\hskip -2em
u_1=\cot\psi,
\mytag{8.6}\\
\hskip -2em
M=\frac{1}{4}\,\bigl(\tan^2\psi\,\tan(2\,\alpha)-4\,\cot(2\,\alpha)\bigr), 
\mytag{8.7}
\endgather
$$
The formula \mythetag{8.7} is derived from the formulas (124) and (1.26) 
in \mycite{4} upon expressing $\alpha_1$ through $\alpha$ by means of 
\mythetag{8.1}. Indeed, we can writhe the code 
\medskip
\parshape 1 20pt 340pt
\noindent
\darkred{{\tt restart:\newline
Eq\_124:=M1\^{}2-M\^{}2-tan(psi)\^{}2:\newline
Eq\_126:=(1-M\^{}2)*sin(2*alpha)-2*M*cos(2*alpha)\newline
\phantom{aaaaaaa}-(1-M1\^{}2)*sin(2*alpha1)+2*M1*cos(2*alpha1):\newline
Eq\_126:=subs(alpha1=Pi/2-alpha,Eq\_126):\newline
sss:=solve(\{Eq\_124,Eq\_126\},\{M,M1\}):\newline
assign(sss):\newline
'M'=M;
}}
\medskip
\noindent
Denoting $s_1=s$ and applying the first formula \mythetag{3.6} with $k=1$,
we get 
$$
\hskip -2em
u_1=\frac{1-s^2}{2\,s}.
\mytag{8.8}
$$
This formula \mythetag{8.8} coincides with the formula \mythetag{6.18}. The 
formulas \mythetag{3.6} with $k=2$, $k=3$, and $k=4$ are written as follows: 
$$
\xalignat 3
&\hskip -2em
u_2=\frac{1-s_2^2}{2\,s_2},
&&u_3=\frac{1-s_3^2}{2\,s_3},
&&u_4=\frac{1-s_4^2}{2\,s_4}.
\quad
\mytag{8.9}
\endxalignat
$$\par
     Now let's recall that $\alpha$ is a Heron angle with the generator 
$m$ (see Definition~\mythedefinition{4.1} and the formula (25) on page 3
of the paper \mycite{4}). Definition~\mythedefinition{4.1} means that the
formula \mythetag{6.12} holds for the angle $\alpha$. From \mythetag{6.12}
we derive the formulas \mythetag{6.13}, \mythetag{6.14}, \mythetag{6.15},
and \mythetag{6.16}. From \mythetag{6.15} and \mythetag{6.16} we derive
$$
\xalignat 2
&\hskip -2em
\tan(2\,\alpha)=
\frac{4\,m\,\bigl(1-m^2\bigr)}
{\bigl(1-m^2\bigr)^2-4\,m^2},
&&\cot(2\,\alpha)=
\frac{\bigl(1-m^2\bigr)^2-4\,m^2}
{4\,m\,\bigl(1-m^2\bigr)}.
\qquad
\mytag{8.10}
\endxalignat
$$
From \mythetag{6.13} and \mythetag{6.14}, applying the formulas \mythetag{4.10},
we derive 
$$
\xalignat 2
&\hskip -2em
\omega_{+}(\alpha)=
\frac{1-m^2+2\,m}{1+m^2},
&&\omega_{-}(\alpha)=
\frac{1-m^2-2\,m}{1+m^2}.
\qquad
\mytag{8.11}
\endxalignat
$$
Finally, substituting \mythetag{8.8} into the formula \mythetag{8.6}, we derive 
$$
\hskip -2em
\tan\psi=\frac{2\,s}{1-s^2}.
\mytag{8.12}
$$
Substituting \mythetag{8.12} and \mythetag{8.10} into \mythetag{8.7} we derive
some definite formula expressing $M$ through $s$ and $m$. This action is performed
by the following code:
\medskip
\parshape 1 20pt 340pt
\noindent
\darkred{{\tt u1:=(1/s-s)/2:\newline
M:=subs(tan(psi)=1/u1,M):\newline
M:=subs(tan(2*alpha)=4*m*(1-m\^{}2)/((1-m\^{}2)\^{}2-4*m\^{}2),M):\newline
M:=subs(cot(2*alpha)=((1-m\^{}2)\^{}2-4*m\^{}2)/4/m/(1-m\^{}2),M):
}}
\medskip
\noindent
Then we use the formulas (140), (141), (142) from \mycite{4}. Applying the formulas
\mythetag{8.11} to them, we derive some definite formulas expressing $u_2$, $u_3$, 
and $u_4$ through $s$ and $m$. This action is performed by means of the following
code:
\medskip
\parshape 1 20pt 340pt
\noindent
\darkred{{\tt u2:=u1*M:\newline
u3:=u1*(omega\_plus(alpha)-M*omega\_minus(alpha)):\newline
u4:=u1*(omega\_minus(alpha)+M*omega\_plus(alpha)):\newline
\newline
u3:=subs(omega\_plus(alpha)=(1-m\^{}2+2*m)/(1+m\^{}2),u3):\newline
u3:=subs(omega\_minus(alpha)=(1-m\^{}2-2*m)/(1+m\^{}2),u3):\newline
u4:=subs(omega\_plus(alpha)=(1-m\^{}2+2*m)/(1+m\^{}2),u4):\newline
u4:=subs(omega\_minus(alpha)=(1-m\^{}2-2*m)/(1+m\^{}2),u4):
}}
\medskip
\noindent
It turns out that the same formulas expressing $u_2$, $u_3$, and $u_4$ 
through $s$ and $m$ can be obtained by substituting \mythetag{6.19},
\mythetag{6.22}, and \mythetag{6.23} into the formulas \mythetag{8.9}. 
This fact confirms once more that the above observation formulated in 
Theorem~\mythetheorem{8.1} is valid. We prove this fact by means of
the following code:
\medskip
\parshape 1 20pt 340pt
\noindent
\darkred{{\tt theta:=(1-s\^{}2)*((1-m\^{}2)\^{}2-4*m\^{}2)/4/m/s/(1-m\^{}2):\newline
eta:=4*m*s*(1-m\^{}2)/(1-s\^{}2)/(1+m\^{}2)/(1-m\^{}2+2*m):\newline
zeta:=(1-s\^{}2)*(1+m\^{}2)*(1-m\^{}2-2*m)/4/m/s/(1-m\^{}2):\newline
\newline
Eq\_u2:=u2-subs(s2=theta,(1-s2\^{}2)/2/s2):\newline
Eq\_u2:=numer(normal(Eq\_u2)):\newline
\newline
Eq\_u3:=u3-subs(s3=eta,(1-s3\^{}2)/2/s3):\newline
Eq\_u3:=numer(normal(Eq\_u3)):\newline
\newline
Eq\_u4:=u4-subs(s4=zeta,(1-s4\^{}2)/2/s4):\newline
Eq\_u4:=numer(normal(Eq\_u4)):\newline
\newline
Eq\_u2,Eq\_u3,Eq\_u4:
}}
\medskip
     Using \mythetag{7.15}, on pages 21 and 22 of his paper \mycite{4} Walter Wyss 
proves that there are no rectangular rational cuboids within his special example
defined by the condition \mythetag{8.1}. Due to Theorem~\mythetheorem{8.1} we see
that the same result is proved in the form of Theorem 2 on pages 12 and 13 of his
paper \mycite{4}.\par
     Walter Wyss's Theorem 2 is valid. It means that there are no rectangular perfect
cuboids within two-dimensional subvarieties $\Gamma_2^1$, $\Gamma_2^2$, $\Gamma_2^3$, 
$\Gamma_2^4$ \pagebreak given by the formulas \mythetag{6.27} and \mythetag{6.31}. 
Neither one of the two-dimensional subvarieties nor their union covers the 
three-dimensional algebraic variety $\Gamma_{3\sssize++}$ given by 
Theorem~\mythetheorem{3.1}. Therefore rectangular perfect cuboids are still possible. 
\head
9. Back to the general case.
\endhead
     On page 22 of his paper \mycite{4} and in Appendix F of this paper Walter Wyss
studies the equation (126). This equation is written as follows:
$$
\hskip -2em
\aligned
(1-M^2)\,\sin(2\,\alpha)&-2\,M\,\cos(2\,\alpha)=\\
&=(1-M_1^2)\,\sin(2\,\alpha_1)-2\,M_1\,\cos(2\,\alpha_1). 
\endaligned
\mytag{9.1}
$$
Denoting through $-4\,D$ the value of each side of the equation \mythetag{9.1},
Walter Wyss splits it into two separate equations:
$$
\gather
(M^2-1)\,\sin(2\,\alpha)+2\,M\,\cos(2\,\alpha)=4\,D,
\mytag{9.2}\\
(M_1^2-1)\,\sin(2\,\alpha_1)+2\,M_1\,\cos(2\,\alpha_1)=4\,D.
\mytag{9.3}
\endgather
$$
The equation \mythetag{9.2} is a quadratic equation with respect to $M$. 
Walter Wyss denotes through $\Delta^2$ the quoter of its discriminant: 
$$
\hskip -2em
\Delta^2=\cos^2(2\,\alpha)+\sin^2(2\,\alpha)+4\,D\,\sin(2\,\alpha).
\mytag{9.4}
$$
The equation \mythetag{9.4} can be derived by means of the following code:
\medskip
\parshape 1 20pt 340pt
\noindent
\darkred{{\tt restart:\newline
Eq\_9\_2:=(M\^{}2-1)*sin(2*alpha)+2*M*cos(2*alpha)-4*D:\newline
Eq\_9\_4:=Delta\^{}2-discrim(Eq\_9\_2,M)/4;\newline
}}
\medskip
\noindent
The equation \mythetag{9.4} is simplified with the use of the well-known
trigonometric identity $\cos^2(2\,\alpha)+\sin^2(2\,\alpha)=1$. As a result we get
$$
\hskip -2em
\Delta^2=1+4\,D\,\sin(2\,\alpha).
\mytag{9.5}
$$
In terms of the machine codes this transformation is performed as follows:
\medskip
\parshape 1 20pt 340pt
\noindent
\darkred{{\tt Eq\_9\_5:=subs(cos(2*alpha)\^{}2=1-sin(2*alpha)\^{}2,Eq\_9\_4);
}}
\medskip
\noindent
The solution of the equation \mythetag{9.2} for $M$ is written as
$$
\hskip -2em
M=\frac{\pm\,\Delta-\cos(2\,\alpha)}{\sin(2\,\alpha)}.
\mytag{9.6}
$$
The formula \mythetag{9.6} is obtained by means of the following code:
\medskip
\parshape 1 20pt 340pt
\noindent
\darkred{{\tt Eq\_9\_2:=subs(D=solve(Eq\_9\_5,D),Eq\_9\_2):\newline
sss:=solve(Eq\_9\_2,M):\newline
M\_plus:=simplify(subs(cos(2*alpha)\^{}2=1-sin(2*alpha)\^{}2,sss[1]))\newline
\phantom{aaaaaaaaa}assuming Delta::positive;\newline
M\_minus:=simplify(subs(cos(2*alpha)\^{}2=1-sin(2*alpha)\^{}2,sss[2]))\newline
\phantom{aaaaaaaaa}assuming Delta::positive;
}}
\medskip
     Then Walter Wyss considers the equation \mythetag{9.5} and writes it as 
follows:
$$
\hskip -2em
(\Delta-1)(\Delta+1)=4\,D\,\sin(2\,\alpha).
\mytag{9.7}
$$
Due to \mythetag{5.9} we know that $\sin(2\,\alpha)\neq 0$. Assume additionally that
$$
\hskip -2em
D\neq 0.
\mytag{9.8}
$$
Under the assumption \mythetag{9.8} we have 
$$
\xalignat 2
&\hskip -2em
\Delta-1\neq 0,
&&\Delta+1\neq 0.
\mytag{9.9}
\endxalignat
$$
Applying \mythetag{9.8} and \mythetag{9.9} to \mythetag{9.7}, we can write it
as follows:
$$
\hskip -2em
\frac{\Delta-1}{4\,D}=\frac{\sin(2\,\alpha)}{\Delta+1}
\mytag{9.10}
$$
The quotients in both sides of \mythetag{9.10} are nonzero. Let's denote their
values through $-r^{-1}/2$, where $r\neq 0$. As a result we split \mythetag{9.10} 
into two equations:
$$
\xalignat 2
&\hskip -2em
\frac{\Delta-1}{4\,D}=-\frac{1}{2\,r},
&\frac{\sin(2\,\alpha)}{\Delta+1}=-\frac{1}{2\,r}.
\mytag{9.11}
\endxalignat
$$
The equations \mythetag{9.11} can be written as linear equations with respect to 
$D$ and $\Delta$:
$$
\xalignat 2
&\hskip -2em
\Delta-1=-\frac{2\,D}{r},
&&\sin(2\,\alpha)=-\frac{\Delta+1}{2\,r}.
\mytag{9.12}
\endxalignat
$$
Resolving the equations \mythetag{9.12}, we get 
$$
\xalignat 2
&\hskip -2em
D=\sin(2\,\alpha)\,r^2+r,
&&\Delta=-2\,\sin(2\,\alpha)\,r-1.
\mytag{9.13}
\endxalignat
$$
The first formula \mythetag{9.13} coincides with the second formula (150) on
page 22 of Walter Wyss's paper \mycite{4}. These two formulas are derived
using the code
\medskip
\parshape 1 20pt 340pt
\noindent
\darkred{{\tt Eq\_9\_11\_1:=(Delta-1)/4/D=-1/2/r;\newline
Eq\_9\_11\_2:=sin(2*alpha)/(Delta+1)=-1/2/r;\newline
sss:=solve(\{Eq\_9\_11\_1,Eq\_9\_11\_2\},\{D,Delta\});
}}
\medskip
Substituting the second formula \mythetag{9.13} into \mythetag{9.6}, we derive 
two solutions for $M$:
$$
\hskip -2em
\aligned
&M_{\sssize+}=-2\,r-\cot(\alpha),\\
\vspace{1ex}
&M_{\sssize-}=2\,r+\tan(\alpha). 
\endaligned
\mytag{9.14}
$$
The formulas \mythetag{9.14} are derived by means of the following code:
\medskip
\parshape 1 20pt 340pt
\noindent
\darkred{{\tt unprotect(D):\newline
assign(sss):\newline
M\_plus:=expand(M\_plus);\newline
M\_minus:=normal(expand(M\_minus)):
}}
\medskip
\parshape 1 20pt 340pt
\noindent
\darkred{{\tt 
M\_minus:=subs(cos(alpha)\^{}2=1-sin(alpha)\^{}2,M\_minus):\newline
M\_minus:=expand(M\_minus);
}}
\medskip
\noindent
The formulas \mythetag{9.14} can be verified by substituting them back to the 
equation \mythetag{9.2} along with the first formula \mythetag{9.13}. This is
done by means of the code 
\medskip
\parshape 1 20pt 340pt
\noindent
\darkred{{\tt simplify(expand(subs(M=M\_plus,Eq\_9\_2)),trig),\newline
simplify(expand(subs(M=M\_minus,Eq\_9\_2)),trig);
}}
\medskip
\noindent
Walter Wyss presents only the second formula \mythetag{9.14} for $M$ on page 22 
of his paper \mycite{4} and actually he does not exploit it.\par
     The equations \mythetag{9.2} and \mythetag{9.3} are similar to each other.
Using this analogy, we can write the following formulas similar to
\mythetag{9.13}:
$$
\xalignat 2
&\hskip -2em
D=\sin(2\,\alpha_1)\,r_1^2+r_1,
&&\Delta_1=-2\,\sin(2\,\alpha_1)\,r_1-1.
\mytag{9.15}
\endxalignat
$$
Though being different, the equations \mythetag{9.2} and \mythetag{9.3} share the
same value of $D$. Therefore from \mythetag{9.13} and \mythetag{9.15} we derive
the equation
$$
\hskip -2em
\sin(2\,\alpha)\,r^2+r=\sin(2\,\alpha_1)\,r_1^2+r_1.
\mytag{9.16}
$$
Factoring both sides of \mythetag{9.16}, we get the equation
$$
\hskip -2em
r\,(r\,\sin(2\,\alpha)+1)=r_1\,(r_1\,\sin(2\,\alpha_1)+1).
\mytag{9.17}
$$
\par
     From the inequalities \mythetag{5.9} we conclude that $\sin(2\,\alpha)\neq 0$
and $\sin(2\,\alpha_1)\neq 0$. Moreover, from the formula \mythetag{9.11} we derive 
$r\neq 0$. Similarly $r_1\neq 0$. Both sides of \mythetag{9.17} are equal to $D$,
where $D\neq 0$ (see \mythetag{9.8}). Hence
$$
\xalignat 2
&\hskip -2em
r\,\sin(2\,\alpha)+1\neq 0,
&&r_1\,\sin(2\,\alpha_1)+1\neq 0.
\mytag{9.18}
\endxalignat
$$
Due to \mythetag{9.18} and the inequalities preceding it, the equation 
\mythetag{9.17} is written as  
$$
\hskip -2em
\frac{r\,\sin(2\,\alpha)+1}{r_1}=\frac{r_1\,\sin(2\,\alpha_1)+1}{r}
\mytag{9.19}
$$
Both sides of \mythetag{9.19} are nonzero. Denoting their value through $1/f$,
we split the equation \mythetag{9.19} into two separate equations:
$$
\xalignat 2
&\hskip -2em
\frac{r_1\,\sin(2\,\alpha_1)+1}{r}=\frac{1}{f},
&&\frac{r\,\sin(2\,\alpha)+1}{r_1}=\frac{1}{f}.
\quad
\mytag{9.20}
\endxalignat
$$
The equations \mythetag{9.20} can be written as two linear equations for $r$ 
and $r_1$:
$$
\xalignat 2
&\hskip -2em
\frac{r}{f}-\sin(2\,\alpha_1)\,r_1=1,
&&\sin(2\,\alpha)\,r-\frac{r_1}{f}=-1.
\quad
\mytag{9.21}
\endxalignat
$$
There are two cases for the equations \mythetag{9.21} --- the regular case and
the singular case. In the regular case we have the inequality
$$
\hskip -2em
\sin(2\,\alpha)\,\sin(2\,\alpha_1)\neq \frac{1}{f^2}.
\mytag{9.22}
$$
In this case the equations \mythetag{9.21} are uniquely solvable. Their solution
is given by 
$$
\hskip -2em
\aligned
&r=\frac{f\,(f\,\sin(2\,\alpha_1)+1)}{1-f^2\,\sin(2\,\alpha_1)\,\sin(2\,\alpha)},\\
\vspace{1ex}
&r_1=\frac{f\,(f\,\sin(2\,\alpha)+1)}{1-f^2\,\sin(2\,\alpha_1)\,\sin(2\,\alpha)}.
\endaligned
\mytag{9.23}
$$
The formulas \mythetag{9.23} can be derived by means of the following code:
\medskip
\parshape 1 20pt 340pt
\noindent
\darkred{{\tt Eq\_9\_20\_1:=(r1*sin(2*alpha1)+1)/r-1/f;\newline
Eq\_9\_20\_2:=(r*sin(2*alpha)+1)/r1-1/f;\newline
sss:=solve(\{Eq\_9\_20\_1,Eq\_9\_20\_2\},\{r,r1\});
}}
\medskip
     The formulas \mythetag{9.23} are written on page 22 of Walter Wyss's paper 
\mycite{4} (see (153) and (154)). They are consistent. Their denominators are nonzero 
due to \mythetag{9.22}. Substituting them back to the equations \mythetag{9.13}, we
derive 
$$
\gather
\hskip -2em
D=\frac{f\,(f\,\sin(2\,\alpha_1)+1)\,(f\,\sin(2\,\alpha)+1)}
{(f^2\,\sin(2\,\alpha_1)\,\sin(2\,\alpha)-1)^2},
\mytag{9.24}\\
\vspace{1ex}
\hskip -2em
\Delta=\frac{f^2\,\sin(2\,\alpha_1)\,\sin(2\,\alpha)+2\,f\,\sin(2\,\alpha)+1}{f^2\,\sin(2\,\alpha_1)\,\sin(2\,\alpha)-1}.
\mytag{9.25}
\endgather
$$
Substituting \mythetag{9.23} into the equations \mythetag{9.15}, we derive 
the same formula \mythetag{9.24} for the parameter $D$ and the following formula 
for $\Delta_1$:
$$
\hskip -2em
\Delta_1=\frac{f^2\,\sin(2\,\alpha_1)\,\sin(2\,\alpha)+2\,f\,\sin(2\,\alpha_1)+1}{f^2\,\sin(2\,\alpha_1)\,\sin(2\,\alpha)-1}.
\mytag{9.26}
$$
The formulas \mythetag{9.24}, \mythetag{9.25}, \mythetag{9.26} are computed
by means of the following code:
\medskip
\parshape 1 20pt 340pt
\noindent
\darkred{{\tt assign(sss):\newline
D:=normal(D);\newline
Delta:=normal(Delta);\newline
\newline
D-normal(sin(2*alpha1)*r1\^{}2+r1);\newline
Delta1:=normal(-2*sin(2*alpha1)*r1-1);
}}
\medskip
     Due to \mythetag{9.14} we have two options for $M$. Similarly, we have two
options for $M_1$:
$$
\hskip -2em
\aligned
&M_{1\sssize+}=-2\,r_1-\cot(\alpha_1),\\
\vspace{1ex}
&M_{1\sssize-}=2\,r_1+\tan(\alpha_1). 
\endaligned
\mytag{9.27}
$$
Substituting \mythetag{9.23} into \mythetag{9.14}, we get the following two
formulas:
$$
\gather
\hskip -2em
M_{\sssize+}=-\frac{2\,f\,(f\,\sin(2\,\alpha_1)+1)}
{1-f^2\,\sin(2\,\alpha_1)\,\sin(2\,\alpha)}-\frac{\cos(\alpha)}{\sin(\alpha)},
\mytag{9.28}\\
\vspace{1ex}
\hskip -2em
M_{\sssize-}=\frac{2\,f\,(f\,\sin(2\,\alpha_1)+1)}
{1-f^2\,\sin(2\,\alpha_1)\,\sin(2\,\alpha)}+\frac{\sin(\alpha)}{\cos(\alpha)}. 
\mytag{9.29}
\endgather
$$
Similarly, substituting \mythetag{9.23} into \mythetag{9.27}, we get the other two
formulas
$$
\gather
\hskip -2em
M_{1\sssize+}=-\frac{2\,f\,(f\,\sin(2\,\alpha)+1)}
{1-f^2\,\sin(2\,\alpha_1)\,\sin(2\,\alpha)}-\frac{\cos(\alpha_1)}{\sin(\alpha_1)},
\mytag{9.30}\\
\vspace{1ex}
\hskip -2em
M_{1\sssize-}=\frac{2\,f\,(f\,\sin(2\,\alpha)+1)}
{1-f^2\,\sin(2\,\alpha_1)\,\sin(2\,\alpha)}+\frac{\sin(\alpha_1)}{\cos(\alpha_1)}. 
\mytag{9.31}
\endgather
$$
The formulas \mythetag{9.28}, \mythetag{9.29}, \mythetag{9.30}, \mythetag{9.31} 
are derived by means of the code
\medskip
\parshape 1 20pt 340pt
\noindent
\darkred{{\tt M\_plus:=M\_plus;\newline
M\_minus:=M\_minus;\newline
M1\_plus:=-cos(alpha1)/sin(alpha1)-2*r1;\newline
M1\_minus:=sin(alpha1)/cos(alpha1)+2*r1;
}}
\medskip
\noindent
The formulas \mythetag{9.28}, \mythetag{9.29}, \mythetag{9.30}, \mythetag{9.31} 
provides four options for choosing the values of $M$ and $M_1$ in the equation 
\mythetag{9.1}: 
$$
\xalignat 2
&\text{1). }M=M_{\sssize+}\text{\ and }M_1=M_{1\sssize+},
&&\text{2). }M=M_{\sssize+}\text{\ and }M_1=M_{1\sssize-},\\
&\text{3). }M=M_{\sssize-}\text{\ and }M_1=M_{1\sssize+},
&&\text{4). }M=M_{\sssize-}\text{\ and }M_1=M_{1\sssize-}.
\endxalignat
$$ 
Choosing any one of these four options, we get a three-parametric solution
of the equation \mythetag{9.1}. {\bf Warning}: a rational solution of the equation 
\mythetag{9.1} does not necessarily produce a rational solution for the cuboid 
equations \mythetag{3.4} and \mythetag{3.5}.\par
     Walter Wyss does not study the singular case for the equations 
\mythetag{9.21} where the inequality \mythetag{9.22} turns to the equality. 
He chooses the option 4. Then on page 22 of his paper \mycite{4}, using the 
formulas \mythetag{9.23}, he considers two cases for $f\neq 0$:
\roster
\rosteritemwd=10pt
\item"(i)" $r_1\neq r$, where he derives the formula
$$
f=\frac{r_1-r}{r\,\sin(2\,\alpha)-r_1\,\sin(2\,\alpha_1)},
\mytag{9.32}
$$
\item"(ii)" $r_1=r$, where he derives the formulas
$$
\gather
f=\frac{r}{1+r\,\sin(2\,\alpha)},
\mytag{9.33}\\
\vspace{1ex}
\sin(2\,\alpha)=\sin(2\,\alpha_1).
\mytag{9.34}
\endgather
$$
\endroster
The formula \mythetag{9.32} is verified by means of the following code:
\medskip
\parshape 1 20pt 340pt
\noindent
\darkred{{\tt Eq\_9\_32:=f-(r1-r)/(r*sin(2*alpha)-r1*sin(2*alpha1)):\newline
Eq\_9\_32:=normal(Eq\_9\_32);
}}
\medskip
\noindent
The formula \mythetag{9.34} is derived from \mythetag{9.23} through
the following equation:
$$
\pagebreak
r_1-r=\frac{f^2\,(\sin(2\,\alpha)-\sin(2\,\alpha_1))}
{1-f^2\,\sin(2\,\alpha_1)\,\sin(2\,\alpha)}=0.
$$
Then the formula \mythetag{9.33} is verified by means of the following code:
\medskip
\parshape 1 20pt 340pt
\noindent
\darkred{{\tt r:=normal(subs(sin(2*alpha1)=sin(2*alpha),r)):\newline
Eq\_9\_33:=f-r/(1+r*sin(2*alpha)):\newline
Eq\_9\_33:=normal(Eq\_9\_33);
}}
\medskip
\noindent
The formula \mythetag{9.34} produces two options:
$$
\xalignat 2
&\hskip -2em
\alpha_1=\alpha,
&&\alpha_1+\alpha=\frac{\pi}{2}.
\mytag{9.35}
\endxalignat
$$
The second option \mythetag{9.35} coincides with the condition \mythetag{8.1}.
\head
10. The rectangular cuboid limit.
\endhead
     The rectangular cuboid limit is the case where the parallelogram ABFE in 
Fig\.~1.1 turns to a rectangle. This case is characterized by the equalities 
$$
\xalignat 2
&\hskip -2em
u_3=u_3,
&&s_3=s_4.
\mytag{10.1}
\endxalignat
$$
The equalities \mythetag{10.1} imply \mythetag{7.17}. Substituting \mythetag{7.17} 
into \mythetag{9.2}, we get
$$
\hskip -2em
D=0. 
\mytag{10.2}
$$
Note that the formulas \mythetag{9.23} were derived under the assumption \mythetag{9.8}.
Comparing \mythetag{9.8} with \mythetag{10.2}, we see that the formulas \mythetag{9.23} 
are not applicable to the rectangular cuboid case directly. For this reasom in his paper \mycite{4} Walter Wyss applies the formulas \mythetag{9.23} through the limit procedure
$$
\xalignat 2
&\hskip -2em
D\neq 0,
&&D\rightarrow 0.
\mytag{10.3}
\endxalignat
$$
Then, using a not well-detailed reasoning, from  $D\rightarrow 0$ in \mythetag{10.3} 
he derives
$$
\xalignat 2
&\hskip -2em
r\rightarrow 0,
&&r_1\rightarrow 0,
\mytag{10.4}\\
\vspace{1ex}
&\hskip -2em
f\rightarrow 0,
&&\frac{r}{r_1}\rightarrow 1.
\mytag{10.5}
\endxalignat
$$
Looking at \mythetag{10.4} and \mythetag{10.5}, he says that $r$ and $r_1$ are
``infinitesimally equal'' and concludes that the rectangular cuboid limit falls
under the case \therosteritem{ii} in \mythetag{9.33} and \mythetag{9.34}. This is 
the crucial mistake in his arguments --- {\bf infinitesimally equal} does not
mean {\bf equal}. There are a lot of rectangular cuboids, no matter rational or
irrational, that do not fall under the case \therosteritem{ii} and \mythetag{9.34}, 
i\.\,e\. such that
$$
\hskip -2em
\sin(2\,\alpha)\neq\sin(2\,\alpha_1). 
\mytag{10.6}
$$
The equality $D=0$ for such cuboids can be reached through the limit procedure 
as $f\rightarrow 0$ in \mythetag{9.23}. Indeed, the formulas \mythetag{9.23} 
simplify to
$$
\pagebreak 
\xalignat 2
&\hskip -2em
r=f+f^2\,\sin(2\,\alpha_1)+o(f^2),
&&r_1=f+f^2\,\sin(2\,\alpha)+o(f^2)
\quad
\mytag{10.7}
\endxalignat
$$
as $f\rightarrow 0$. Due to \mythetag{10.6} $r$ and $r_1$ in \mythetag{10.7}
tend to zero never being equal for sufficiently small $f\neq 0$. Rational 
rectangular cuboids are not yet found.\par
\head
11. Conclusions.
\endhead
     Almost all formulas in Walter Wyss's paper \mycite{4} have been verified.
They are correct. For the reader's convenience all of the code used for verifying
formulas is collected in ancillary files in the section-by-section form according
to the sections of the present paper.\par
     Walter Wyss's paper \mycite{4} comprises a valuable result for the theory of 
slanted cuboids. This result is expressed by the explicit formulas \mythetag{6.24}, 
\mythetag{6.25}, and \mythetag{6.26} that produce four explicit two-parametric 
solutions of the basic slanted cuboid equation \mythetag{3.9} through the formulas
\mythetag{6.27} and \mythetag{6.31}.\par
     As for the main goal of the paper \mycite{4}, it is not reached. The paper 
does not contain a correct proof for the no perfect cuboid claim in its title. 
\par 
\adjustfootnotemark{-1}

\enddocument
\end